\documentclass[12 pt]{amsart}
\usepackage{amsthm, amsfonts, amssymb, color}
\usepackage{mathrsfs}
\usepackage{amsmath}
\usepackage{amstext, amsxtra}
\usepackage{txfonts}
\usepackage[colorlinks, linkcolor=black, citecolor=blue, hypertexnames=false]{hyperref}
\usepackage{comment}

\allowdisplaybreaks
\usepackage{geometry}
\newtheorem{theorem}{Theorem}[section]

\newtheorem{corollary}[theorem]{Corollary}
\newtheorem{lemma}[theorem]{Lemma}
\newtheorem{definition}[theorem]{Definition}
\newtheorem{example}[theorem]{Example}
\newtheorem{assumption}{Assumption}[section]
\newtheorem{remark}[theorem]{Remark}

\newcommand{\jap}[1]{\langle #1\rangle}

\newcommand{\A}{\mathcal{A}}

\newcommand{\eq}{\begin{equation}}
\newcommand{\eeq}{\end{equation}}
\newcommand{\eqa}{\begin{eqnarray}}
\newcommand{\eeqa}{\end{eqnarray}}

\newcommand\Hi{\mathcal{H}}

\newcommand\R{\mathbb{R}}

\newcommand\vu{\vec{u}}

\numberwithin{equation}{section}

\title   [] { On The large Time Asymptotics of  Klein-Gordon type equations with General Data}

\author{   Avy Soffer }

\address{Department of Mathematics\\
Rutgers University\\
110 Frelinghuysen Rd.\\
Piscataway, NJ, 08854, USA}

\email{soffer@math.rutgers.edu}

\author{Xiaoxu Wu}
\address{Mathematical Sciences Institute\\
Australian National University\\
Acton, ACT 2601, Australia}
\email{Xiaoxu.Wu@anu.edu.au}

\thanks{
A.Soffer is supported in part by NSF Grant DMS-2205931
}

\begin{document}

\begin{abstract}
We study the Klein-Gordon equation with general interaction terms, which may be linear or nonlinear,
 and space-time dependent. We initiate the study of such equations with large (non-radial) data.   We prove that global solutions are asymptotically given by a free wave and a  weakly localized part.
The proof is based on constructing in a new way the Free Channel Wave Operator, and further tools from the recent works \cite{Liu-Sof1,Liu-Sof2,SW2020,SW2022}.
This work generalizes the results of part of  \cite{Liu-Sof1,Liu-Sof2} on the Schr\"odinger equation to arbitrary dimension, and non-radial data.
\end {abstract}
\maketitle

\section{Introduction}

\subsection{Setup. }We start by describing the setup and main results of this paper, before moving to the background and discussions in Section \ref{sec: back}. Let $\langle \cdot\rangle:\mathbb{R}^n\to\mathbb{R}$, $x\mapsto\sqrt{1+|x|^2}$. In space dimension $n\geq1$, we consider the following nonlinear Klein-Gordon (KG) equations:
\eq\tag{KG}
\begin{cases}
(\square+1)u=N(u,x,t)\,u\\
\vu(0):=(u(0),\dot{u}(0))=(u_0, \dot{u}_0)\in \Hi
\end{cases}, \quad (x,t)\in \mathbb{R}^n\times \mathbb{R}, \label{KG}
\eeq
where $p:=-i\nabla_x$, so that $\jap{p}:=\sqrt{1+|p|^2}=\sqrt{H_0+1}$ with $H_0:=-\Delta_x$; we use the two notations $\jap{p}$ and $\sqrt{H_0+1}$ interchangeably throughout. The space $\Hi:=H^1(\mathbb{R}^n)\times L^2(\mathbb{R}^n)$ denotes a Hilbert space equipped with the inner product
\eq\label{eq: inner product}
(\vec{u}, \vec{v})_{\Hi}:=(\jap{p} u_1, \jap{p} v_1 )_{L^2_x(\mathbb{R}^n)}+( u_2,  v_2 )_{L^2_x(\mathbb{R}^n)}\quad \text{ for all }\vec{v}, \vec{u}\in \Hi.
\eeq
The inner product \eqref{eq: inner product} induces the norm $\|\vu\|_{\Hi}^2:=(\vu,\vu)_{\Hi}$. Throughout, we write $u(t):=u(\cdot,t)$ for the solution at time $t$ regarded as a function of $x$, $\dot u(t):=\partial_t u(\cdot,t)$ for its time derivative, and $\vec u(t):=(u(t),\dot u(t))\in\Hi$ for the pair. Here $\Box:=\partial^2_t-\Delta_x$ and $N$ represents the interaction, which will be detailed later. We write $X\lesssim Y, Y\gtrsim X$ to indicate $X\leq CY$ for some constant $C>0$ and $X\lesssim_a Y$ to indicate $X\leq CY$ for some $C=C(a)>0$. Unsubscripted norms always denote the $L^2_x(\mathbb{R}^n)$ norm or, when applied to an operator, the $L^2_x(\mathbb{R}^n)\to L^2_x(\mathbb{R}^n)$ operator norm; that is, $\|\cdot\|:=\|\cdot\|$ on functions and $\|\cdot\|:=\|\cdot\|$ on operators. Throughout, $c_n:=(2\pi)^{-n/2}$, and for $f,g\in \mathcal{S}(\mathbb{R}^n)$ (the Schwartz class) the Fourier transform and its inverse are taken with the symmetric convention
\eq\label{eq: Fourier}
\hat f(\xi):=c_n\!\int_{\mathbb{R}^n}e^{-ix\cdot\xi}f(x)\,dx,\qquad \check g(x):=c_n\!\int_{\mathbb{R}^n}e^{ix\cdot\xi}g(\xi)\,d\xi,
\eeq
extended by density to tempered distributions in the standard way.

In this paper, we will study the long-time behavior of global solutions to \eqref{KG} as $t\to \infty$. We aim to demonstrate a decomposition result for these global solutions.
\subsection{Assumptions and examples}We consider initial data $\vec{u}(0)$ that lead to a global solution in $\Hi$:
\begin{assumption}\label{asp: global}We assume that $\vec{u}(0)\in \Hi$ leads to a global solution in $\Hi$ whose $\Hi$-norm is uniformly bounded in time:
\begin{equation*}
E:=\sup\limits_{t\geq 0} \| \vec{u}(t)\|_{\Hi}<\infty.
\end{equation*}
\end{assumption}
We consider two types of $N(u,x,t)$:
\begin{assumption}\label{asp: 2}(Theorem \ref{thm1}, \,\, $n\geq 1$)  {\bf Local interactions}. For some $\sigma>1$, we assume $\langle x\rangle^\sigma N(u(t),x,t) u(t)$ remains bounded in $L^2_x(\mathbb{R}^n)$ uniformly in time:
\begin{equation*}
E_u:=\sup\limits_{t\geq 0}\| \langle x\rangle^\sigma N(u(t),x,t)u(t) \|<\infty.
\end{equation*}

\end{assumption}
\begin{assumption}\label{asp: 3}(Theorem \ref{thm2}, $n\geq 3$) {\bf Non-local interactions}. $N(u(t),x,t)u(t)$ remains bounded in $L^1_x(\mathbb{R}^n)$ uniformly in time:
\begin{equation*}
\sup\limits_{t\geq 0}\| N(u(t),x,t)u(t)\|_{L^1_x(\mathbb{R}^n)}<\infty.
\end{equation*}
\end{assumption}
When $N$ depends on $u$, Assumption~\ref{asp: global} fixes the orbit $\{u(t)\}_{t\geq 0}$ along which the equation is studied; we may then regard $N(u(t),x,t)$ as a time-dependent multiplication operator along this orbit.
\begin{example}[Local interactions]\label{ex: local}
A typical local interaction in space dimension $n=1$ is
\begin{equation*}
N(u,x,t)= V(x,t)+a(x)u+b(x)u^2,
\end{equation*}
where $\langle x\rangle^\sigma V(x,t)\in L^\infty_{x,t}(\mathbb{R}\times\mathbb{R}_{\geq 0})$ and $\langle x\rangle^\sigma a(x), \langle x\rangle^\sigma b(x)\in L^\infty_x(\mathbb{R})$ for some $\sigma>1$.
\end{example}
\begin{example}[Non-local interactions]
A typical non-local interaction in space dimension $n\geq 3$ is
\begin{equation*}
N(u,x,t)= V(x,t)+\lambda u^2+\lambda' u^3,
\end{equation*}
where $V(x,t)\in L^\infty_tL_x^2(\mathbb{R}^n\times\mathbb{R}_{\geq 0})$ and $\lambda,\lambda'\in \mathbb{R}$.
\end{example}
\begin{example}[Charge transfer interactions]\label{ex: charge}More generally, in dimension $n\geq 3$, one can take
\begin{equation*}
N(u,x,t)=V(x,t)+\sum\limits_{j=1}^{K}\epsilon_j \lambda_j |u|^{p_j},\quad \lambda_j>0,\ \epsilon_j\in\{+1,-1\},\ 1\leq p_j\leq \tfrac{n+2}{n-2},
\end{equation*}
where $V(x,t)$ is of charge transfer type, namely
\begin{equation*}
V(x,t)=\sum\limits_{j=1}^M V_j(x-g_j(t)v_j,t),\quad V_j\in L^\infty_t L^2_x,\ v_j\in \mathbb{R}^n,
\end{equation*}
with real-valued $g_j(t)$. See Corollary~\ref{app3} for more details.
\end{example}
\subsection{Free channel wave operators}\label{subsecL free channel}In this section we introduce the key notion for separating the free solution from $\vec u(t)$: the free channel wave operators. Following the notation of \cite{SW2022}, we first introduce the smooth cut-offs $F_c$ and $F_1$.
\begin{definition}\label{def: cutoffs}We denote by $F_c$ and $F_1$ smooth cut-off functions satisfying
\begin{equation*}
F_j(\lambda)=\begin{cases}1 & \text{when }\lambda\geq 1\\ 0 & \text{when }\lambda<1/2\end{cases},\qquad j=c,1.
\end{equation*}
Moreover, for all $a>0$ and $j=c,1$ we define
\begin{equation*}
F_j(\tfrac{\lambda}{a}\leq 1)\equiv 1-F_j(\tfrac{\lambda}{a}),\qquad F_j(\tfrac{\lambda}{a}>1)\equiv F_j(\tfrac{\lambda}{a}).
\end{equation*}
\end{definition}

We now define a free channel wave operator~$\Omega_{\alpha}^*$ acting on initial data $\vec u(0)$ as follows:
 \begin{equation}\label{eq: subsecL free channel, eq 1}
 \Omega_{\alpha}^*\vu(0) := s\text{-}\lim_{t\to \infty} \mathcal{F}_\alpha(x,p,t) U_0(0,t) \vu(t)\quad \text{ in }\Hi,
 \end{equation}
where $U_0(t,0)$ denotes the free solution operator of system \eqref{KG} and $\mathcal{F}_\alpha$ is the matrix operator
 \begin{equation*}
\mathcal{F}_\alpha(x,p,t)=\begin{pmatrix}
 \frac{1}{\jap{p}} F_\alpha(x,p,t)\jap{p} &0 \\
 0 &  F_\alpha(x,p,t)
 \end{pmatrix}.
 \end{equation*}
The function $F_\alpha$ is built from the cut-offs of Definition~\ref{def: cutoffs} and satisfies:
 \begin{enumerate}
\item When $N$ is a localized interaction and $n\geq 1$, $F_\alpha$ is defined as
\eq\label{Falpha1}
     F_\alpha(x,p,t):=  F_c(\frac{|x|}{t^\alpha}\leq 1)F_1(t^{\beta}|p|> 1)F_1(|p|\leq t^{\beta})\quad  
 \eeq
 for $\alpha\in (0,1)$ and $\beta\in (0, \min\{1-\alpha, \frac{\sigma-1}{\sigma}\})$.
 \item When $N$ is a non-local interaction and $n\geq 3$, $F_\alpha$ is defined as
 \eq\label{Falpha2}
     F_\alpha(x,p,t):= F_c(\frac{|x|}{t^\alpha}\leq 1)F_1(|p|\leq t^{\beta})  
 \eeq
 for $\alpha\in (0,1-2/n)$ and $\beta\in (0,\frac{n}{n+3}(1-\alpha)-\frac{2}{n+3})$.
 \end{enumerate}
 \begin{remark}The choice of $\alpha$ and $\beta$ here comes from the dispersive estimates for the free flow, which are proved in Lemma \ref{lem: free decay}:
 \begin{itemize}
     \item When $\alpha\in (0,1)$ and $\beta\in (0, \min\{1-\alpha, \frac{\sigma-1}{\sigma}\})$, the weighted $L^2$ estimate holds:
 \begin{equation*}
 \| F_\alpha(x,p,t)e^{\pm it\jap{p}}\langle x \rangle^{-\sigma} \|\in L^1_t[1,\infty).
 \end{equation*}
  \item  When $\alpha \in (0,1-2/n)$ and $\beta\in (0, \frac{n}{n+3}(1-\alpha)-\frac{2}{n+3})$, the $L^1$ estimate holds:
  \begin{equation*}
  \| F_\alpha(x,p,t) e^{\pm it\jap{p}}\|_{L^1_x(\mathbb{R}^n)\to L^2_x(\mathbb{R}^n)}\in L^1_t[1,\infty).
  \end{equation*}
 \end{itemize}
 \end{remark}
The existence of $\Omega_\alpha^*\vec{u}(0)$ in $\Hi$ is proved in Theorems~\ref{thm1} and~\ref{thm2}.

 \subsection{Main results}
Here are our main results:

\begin{theorem}[Local interactions]\label{thm1}
Let $\vu(t)=(u(t),\dot u(t))$ be a solution to \eqref{KG} satisfying Assumptions~\ref{asp: global} and~\ref{asp: 2}. Denote by $\mathscr{I}_{\mathrm{loc}}$ the set of admissible parameters
\begin{equation*}
\mathscr{I}_{\mathrm{loc}}:=\{(\alpha,\beta)\in (0,1)\times(0,1)\,:\, \beta<\min\{1-\alpha,\tfrac{\sigma-1}{\sigma}\}\}.
\end{equation*}
Then for every $(\alpha,\beta)\in\mathscr{I}_{\mathrm{loc}}$, the free channel wave operator acting on initial data $\vec u(0)$, $\Omega_\alpha^*\vec u(0)$, defined in \eqref{eq: subsecL free channel, eq 1}, exists in $\Hi$, and is independent of $(\alpha,\beta)$:
\begin{equation*}
\Omega_\alpha^*\vec u(0)=\Omega_{\alpha'}^*\vec u(0)\quad\text{for all }(\alpha,\beta),(\alpha',\beta')\in\mathscr{I}_{\mathrm{loc}}.
\end{equation*}
\end{theorem}
\begin{theorem}[Non-local interactions]\label{thm2}
Let $\vu(t)=(u(t),\dot u(t))$ be a solution to \eqref{KG} satisfying Assumptions~\ref{asp: global} and~\ref{asp: 3}. Denote by $\mathscr{I}_{\mathrm{nl}}$ the set of admissible parameters
\begin{equation*}
\mathscr{I}_{\mathrm{nl}}:=\{(\alpha,\beta)\in (0,1-2/n)\times(0,1)\,:\, \beta<\tfrac{n}{n+3}(1-\alpha)-\tfrac{2}{n+3}\}.
\end{equation*}
Then for every $(\alpha,\beta)\in\mathscr{I}_{\mathrm{nl}}$, the free channel wave operator acting on initial data $\vec u(0)$, $\Omega_\alpha^*\vec u(0)$, defined in \eqref{eq: subsecL free channel, eq 1}, exists in $\Hi$, and is independent of $(\alpha,\beta)$:
\begin{equation*}
\Omega_\alpha^*\vec u(0)=\Omega_{\alpha'}^*\vec u(0)\quad\text{for all }(\alpha,\beta),(\alpha',\beta')\in\mathscr{I}_{\mathrm{nl}}.
\end{equation*}
\end{theorem}
After removing all free waves from the solution $\vu(t)$, we further characterize the properties of the remaining component. We will show that if $N(u,x,t)$ is local in space and satisfies Assumption~\ref{asp: 2}, then the expected value of $\langle x\rangle$ with respect to the non-free component cannot asymptotically grow faster than $c\min\{\sqrt t,\,t^{1/\sigma}\}$ for some $c>0$; see \eqref{conlusion: wlocal} below for the precise statement. Let $\mathcal{X}$ denote a matrix operator:
\begin{equation*}
\mathcal{X}:=\begin{pmatrix}
   \jap{p}^{-1} \langle x\rangle \jap{p} & 0\\
   0 & \langle x\rangle
\end{pmatrix}.
\end{equation*}
\begin{theorem}\label{thm3}Assume all assumptions in Theorem \ref{thm1} are satisfied. The solution $\vu(t)$ has the following asymptotic decomposition: as $t\to \infty$, 
\eq\label{limit}
\| \vu(t)-U_0(t,0)\Omega_\alpha^*\vec{u}(0)-\vu_{wlc}(t)\|_{\Hi}\to 0
\eeq 
where $\vec{u}_{wlc}(t)$ is a sub-ballistic state in the sense that, for every $\epsilon\in (0,1/2)$,
\eq\label{conlusion: wlocal}
(\vec{u}_{wlc}(t), \mathcal{X}^N \vec{u}_{wlc}(t))_\Hi \lesssim_{E,\epsilon} \max\bigl\{t^{N(1/2+\epsilon)},\, t^{N(1/\sigma+\epsilon)}\bigr\}
\eeq
for all $N\geq 0$.
\end{theorem}

\subsection{The background.}\label{sec: back}Theorems~\ref{thm1}--\ref{thm3} above fit into a broader research program on the long-time asymptotics of dispersive and hyperbolic wave equations, central to the analysis of evolution equations in Physics and Geometry. It is well known that the asymptotic solutions of such equations, when they exist, exhibit a dizzying zoo of possible structures: besides the \emph{free wave} (solution of the equation without interaction), a multitude of other solutions may appear, localized around possibly moving centers of mass --- nonlinear bound states, solitons, breathers, hedgehogs, vortices, etc.\ \cite{Sof} --- and their analysis is usually case-by-case. A natural question follows: do solutions of dispersive/hyperbolic equations converge, in an appropriate norm ($L^2$ or $H^1$), to a free wave plus independently moving localized parts? This is precisely the statement of \emph{asymptotic completeness}, established in $N$-body scattering \cite{SS1987,Gr2,D-Ger,Dere,SSjams,SSDuke, H-Sig1,H-Sig2}, where the possible outgoing clusters are identified as bound states of subsystems, and more recently for nonlinear completely integrable equations in one dimension \cite{SG,dNLS}. But when the interaction includes time-dependent potentials (even localized in space) or more general nonlinear terms, no a-priori knowledge of the asymptotic states is available.

For time-independent interactions, spectral theory resolves the question: scattering states evolve from the continuous spectrum and the localized part is formed by the point spectrum. Once the interaction is time-dependent or nonlinear, this decomposition is unavailable. Recent works on Schr\"odinger type problems \cite{Liu-Sof1,Liu-Sof2,SW2022} have obtained general results in this setting; the present work initiates the corresponding study for hyperbolic equations, focusing on the KG equation in arbitrary dimension with general interactions, including semi-linear ones. A parallel development for the nonlinear wave equation, sharing several of the propagation-wave arguments adopted here, has been carried out in \cite{MSW2025}. To our knowledge, Theorems~\ref{thm1}--\ref{thm3} provide the first large-data multichannel scattering results for nonlinear KG equations outside the integrable class.

For earlier works on time-dependent potentials we mention: charge transfer Hamiltonians \cite{Yaj1, Gr1,Wul,Per, RSS, RSS2,DSY}, decaying-in-time and small potentials \cite{How, RS}, time periodic potentials \cite{Yaj2, How}, and random (in time) potentials \cite{BeSof}; see also \cite{BeSof1}. For potentials with asymptotic energy distribution more could be done \cite{SS1988}. A recent progress for more general localized potentials without smallness assumptions is obtained in \cite{SW2020}, some tools from which will be used here. Turning to the nonlinear case, Tao \cite{T1,T2,T3} has shown that the asymptotic decomposition holds for NLS with inter-critical nonlinearities, in 3 or higher dimensions, with radial initial data. In particular, in very high dimension and with an interaction that is a sum of a smooth compactly supported potential and a repulsive nonlinearity, Tao showed that the localized part is smooth and localized; in other cases, the localized part is only weakly localized and smooth. Tao's work uses direct estimates of the incoming and outgoing parts of the solution to control the nonlinear term via Duhamel representation, in a certain sense in the spirit of Enss work; see also \cite{Rod-Tao}. For the critical power wave equations and wave map problems there has been great progress in understanding the large time behavior with large data; see e.g.~\cite{DKM2,DKMM,NLWpotential,cote} and cited references. For the nonlinear KG equation, there are no results on multichannel scattering with large data. There are also many major works on the stability of coherent structures, e.g.~\cite{LuSc,GerPu3,Mun,SofWei1}, and a large literature on NLS, KdV and more; for small data and long-range type interactions, see \cite{GerPu1,GerPu2,Lind-S1,Lind-S2,Lind-S3,IT2015,Lind-J-S,Lind-J-S1}.

In contrast, the new approach of Liu-Soffer \cite{Liu-Sof1,Liu-Sof2} is based on proving a-priori estimates on the full dynamics that hold in suitably defined domains of the extended phase-space, namely propagation estimates in domains exterior to the support of the interaction. Similar propagation observables were used in many other (mostly linear) problems with time independent potentials; see e.g.~\cite{HWSS1999,Gr2,D-Ger} and cited references. In this way it was possible to show the asymptotic decomposition for general localized interactions, including time- and space-dependent ones, under radial initial data, ensuring the localization of the nonlinear part. More detailed information can be obtained on the localized part of the solution: besides being smooth, its expanding part (if it exists) can grow at most like $|x|\leq \sqrt t$, and is concentrated in a thin set of the extended phase-space. The free part concentrates on the \emph{propagation set} where $x=vt$, $v=2p$, with $p=-i\nabla_x$ the momentum, while the weakly localized part is localized in the regions
$$|x|\sim t^\alpha \quad\textbf{and} \quad |p|\sim  t^{-\alpha}, \quad\forall ~0\leq\alpha \leq 1/2,$$
showing that the spreading part follows a self-similar pattern; see \cite{Liu-Sof1,Liu-Sof2}. The method of proof consists of three steps: construct the Free Channel Wave Operator; prove localization of the remainder localized part and use it to prove smoothness; and via further propagation estimates adapted to localized solutions, prove concentration on thin sets corresponding to self-similar solutions. In this work we will mainly do the first part and some of the second of \cite{Liu-Sof1,Liu-Sof2}. We emphasize that the spreading localized solutions, if they exist, have a non-small nuclei part around the origin (in both Tao's and Liu-Soffer's results), so these are not pure self-similar solutions, as they appear in the special cases of critical nonlinearities; see e.g.~\cite{S2021,DKM2}. We expect similar behavior for the weakly localized part of solutions to KG equations.

This point of view was generalized in \cite{SW2022} to include non-radial data and interactions, and with localized interactions to arbitrary dimension, in the case of the Schr\"odinger problem. This generalization is based on refined localization of the channel wave operators: by localizing around the phase-space support of the free wave, one obtains a sharper decomposition of the localized and scattering parts, avoiding the need for localization of the interactions in some cases. The idea of sharp localization was used in other ways in long-range scattering theory, e.g.~in \cite{Sig,SSInvention}. Here we follow these constructions for the KG case, mainly by viewing the dynamics of the KG equation as generated by a couple of Schr\"odinger type equations, with dispersion relation $\sqrt{p^2+1}$ and group velocity $v=p/\sqrt{p^2+1}$. The extension of the previous methods proceeds along similar ideas, at least when the interaction terms are localized in space. When the interaction terms are not localized but only satisfy $L^p$ decay conditions (with $1\leq p\leq 2$), the situation is more complicated: $L^p$ decay estimates for the KG equation require control of derivatives of the initial data, due to the poor dispersion for hyperbolic equations near infinite frequency. We deal with this by introducing an extra cutoff of high frequencies into the construction of the channel wave operators, as explained next.

We close this background section by placing our construction in \S\ref{subsecL free channel} within the abstract framework of scattering theory. The key object in multichannel scattering is the \emph{channel wave operator}
\begin{equation}\label{abstract: chan op}
\Omega_a^*\equiv s\text{-}\lim\limits_{t \to \infty} U_a(-t)U(t)\vec{u}(0),
\end{equation}
where the limit is taken in the strong sense in a suitable Hilbert space, $U(t)\vec u(0)$ is the global solution of the KG equation with initial data $\vec u(0)=(u(0),\dot u(0))$ and full dynamics $U(t)=U(t,0)$ generated by the Hamiltonian $H=-\Delta+1+N(u,x,t)$, and $U_a$ is the asymptotic dynamics generated by a channel Hamiltonian $H_a$ associated with a given channel $a$. The free channel corresponds to $H_a=-\Delta+1$, which is the only channel constructed in this work. A crucial observation \cite{SS1987} is that one can equivalently modify \eqref{abstract: chan op} by inserting an auxiliary localization $J_a$:
\begin{equation}\label{abstract: chan op2}
\Omega_a^*\equiv s\text{-}\lim\limits_{t \to \infty} U_a(-t)J_a U(t)\vec{u}(0).
\end{equation}
Our definition \eqref{eq: subsecL free channel, eq 1} of $\Omega_\alpha^*\vec u(0)$ is precisely an instance of \eqref{abstract: chan op2}, with $J_a$ replaced by the phase-space cut-off $\mathcal{F}_\alpha(x,p,t)$ from \eqref{eq: subsecL free channel, eq 1}-\eqref{Falpha2}; the high-frequency cut-off contained in $\mathcal{F}_\alpha$ is precisely what compensates for the poor dispersion of the KG flow at infinite frequency mentioned above.

\section{Preliminaries}
\subsection{Free and perturbed KG: dispersive estimates and Duhamel formula}

Let $\vu_0(t):=(u_0(t),\dot u_0(t))$ solve the free KG equation
\eq\label{Pfree}
\begin{cases}
(\square+1)u_0(t)=0,\\
\vu_0(0)=\vu(0)=(u(0),\dot u(0))\in \Hi,
\end{cases}\qquad (x,t)\in \mathbb{R}^n\times \mathbb{R}.
\eeq
Setting $\A_0:=\begin{pmatrix}0 & -1\\ H_0+1 & 0\end{pmatrix}$, the system \eqref{Pfree} is equivalent to $\partial_t\vu_0(t)=-\A_0\vu_0(t)$, so
\eq\label{Ufree}
\vu_0(t)=e^{-t\A_0}\vu(0)=:U_0(t,0)\vu(0),
\eeq
and the components admit the standard representation
\eq\label{free:rep1}
u_0(t)= \cos(t\jap{p})u(0)+\tfrac{\sin(t\jap{p})}{\jap{p}}\dot u(0),
\eeq
\eq\label{free:rep2}
\dot u_0(t)=-\sin(t\jap{p})\jap{p}\,u(0)+\cos(t\jap{p})\dot u(0).
\eeq
Our analysis relies on the following standard dispersive estimate for the KG propagator (see H\"ormander \cite[Corollary 7.2.4]{Horm}).
\begin{lemma}\label{dlem1}
For all $t\in \mathbb{R}$ and Schwartz $f$ on $\mathbb{R}^n$,
\begin{equation*}
\| e^{\pm it\jap{p}}f\|_{L^\infty_x(\mathbb{R}^n)}\lesssim \tfrac{1}{\langle t\rangle^{n/2}}\,\|\jap{p}^{(n+3)/2} f\|_{L^1_x(\mathbb{R}^n)}.
\end{equation*}
\end{lemma}
We also need the following commutator estimate from \cite[Lemma~3.6]{SW2022}.
\begin{lemma}[{\cite[Lemma~3.6]{SW2022}}]\label{com}Let $F_1^{(l)}(k):=\frac{d^l F_1}{dk^l}$, $l=0,1$. For $t\geq 1$, $\alpha>0$ and $\alpha>\beta$, the following inequality holds for all $n\geq 1$:
\begin{equation*}
\|[F_c(\tfrac{|x|}{t^\alpha}\leq 1),F_1^{(l)}(t^\beta|p|>1)]\|\lesssim_n \tfrac{1}{t^{\alpha-\beta}}.
\end{equation*}
\end{lemma}

\begin{lemma}[Dispersive estimates for free flows]\label{lem: free decay}Let $F_\alpha$ be as in \eqref{Falpha1}, \eqref{Falpha2}, depending on whether $\mathcal{N}$ is a local or non-local interaction, respectively.
\begin{enumerate}
\item (\emph{Local case}) When $\alpha\in (0,1)$ and $\beta\in (0,\min\{1-\alpha,\frac{\sigma-1}{\sigma}\})$,
\eq\label{decay: eq1}
\|F_\alpha(x,p,t)\,e^{\pm it\jap{p}}\,\langle x\rangle^{-\sigma}\|\in L^1_t[1,\infty).
\eeq
\item (\emph{Non-local case}) When $\alpha\in (0,1-2/n)$ and $\beta\in (0,\frac{n}{n+3}(1-\alpha)-\frac{2}{n+3})$,
\eq\label{decay: eq2}
\|F_\alpha(x,p,t)\,e^{\pm it\jap{p}}\|_{L^1_x(\mathbb{R}^n)\to L^2_x(\mathbb{R}^n)}\in L^1_t[1,\infty).
\eeq
\end{enumerate}
\end{lemma}
\begin{proof}\noindent\fbox{\emph{Local case.}} Writing $F_1=F_1(t^\beta|k|>1)F_1(|k|\leq t^\beta)$ for the frequency factor in $F_\alpha$, the phase $\pm t\sqrt{|k|^2+1}$ satisfies $|\nabla_k[t\sqrt{|k|^2+1}]|\geq C_1 t^{1-\beta}> C_1 t^\alpha$ on the support of $F_1$ (since $\beta<1-\alpha$). Split $F_\alpha e^{\pm it\jap{p}}\langle x\rangle^{-\sigma}=Q_1+Q_2$ where $Q_1:=F_\alpha e^{\pm it\jap{p}}\langle x\rangle^{-\sigma}\chi(|x|>\tfrac{C_1}{100}t^{1-\beta})$ (cutoff multiplied on the \emph{right}, i.e., after $\langle x\rangle^{-\sigma}$) and $Q_2:=F_\alpha e^{\pm it\jap{p}}\langle x\rangle^{-\sigma}\chi(|x|\leq\tfrac{C_1}{100}t^{1-\beta})$. The weight gives $\|Q_1\|\lesssim t^{-\sigma(1-\beta)}\in L^1_t[1,\infty)$ since $\sigma(1-\beta)>1$, and non-stationary phase (integration by parts in $k$ inside $Q_2$) gives $\|Q_2\|\lesssim_M t^{-(1-\beta)M}\in L^1_t$ for $M$ large, proving \eqref{decay: eq1}.

\noindent\fbox{\emph{Non-local case.}} By Lemma~\ref{dlem1},
\begin{equation*}
\|F_\alpha(x,p,t)\,e^{\pm it\jap{p}}\|_{L^1_x\to L^2_x}\lesssim t^{\alpha n/2}\cdot\langle t\rangle^{-n/2}\cdot\langle t\rangle^{(n+3)\beta/2},
\end{equation*}
and the assumed range of $\beta$ gives $\tfrac{n}{2}-\tfrac{n\alpha}{2}-\tfrac{n+3}{2}\beta>1$, proving \eqref{decay: eq2}.
\end{proof}
A quick application of Lemma~\ref{lem: free decay} is the dispersive estimate for the matrix operator of free flows. Set
\eq\label{def: matrixN}
\mathcal{N}(\vec{u}(t)):=\begin{pmatrix} 0&0\\ N(u(t),x,t)& 0\end{pmatrix},
\eeq
and let
\eq
B_1(t):=\begin{pmatrix} \jap{p}^{-1}F_1F_\alpha(x,p,t)\jap{p}  & 0\\
0& F_1F_\alpha(x,p,t)
\end{pmatrix}\label{B1}
\eeq
and 
\eq
B_2(t):=\begin{pmatrix}
   \jap{p}^{-1} F_\alpha(x,p,t)F_c(\frac{|x|}{t^\alpha }\leq 1)\jap{p}  & 0\\
   0 & F_\alpha(x,p,t)F_c(\frac{|x|}{t^\alpha }\leq 1)
\end{pmatrix},\label{B2}
\eeq
where $F_1$ is given by
\begin{equation*}
F_1=\begin{cases}
    F_1(t^\beta|p|>1)F_1(|p|\leq t^\beta)& \text{ if }N\text{ is local} \\
    F_1(|p|\leq t^\beta)& \text{ if }N\text{ is non-local}
\end{cases}.
\end{equation*}
\begin{corollary}[Dispersive estimates for free flows]\label{cor: dispersive}Suppose Assumption~\ref{asp: global} holds. Then for $j=1,2$,
\eq\label{Cor: state1}
\|B_j(t) U_0(0,t)\mathcal{N}(\vec{u}(t))\vec{u}(t)\|_{\mathcal{H}}\in L^1_t[1,\infty)
\eeq
and
\eq\label{Cor: state2}
\|\mathcal{F}_\alpha(x,p,t) U_0(0,t)\mathcal{N}(\vec{u}(t))\vec{u}(t)\|_{\mathcal{H}}\in L^1_t[1,\infty),
\eeq
provided either
\begin{enumerate}
\item[(i)] (\emph{Local case}) Assumption~\ref{asp: 2} holds, or
\item[(ii)] (\emph{Non-local case}) Assumption~\ref{asp: 3} holds and $n\geq 3$.
\end{enumerate}
\end{corollary}
\begin{proof}By \eqref{free:rep1} and \eqref{free:rep2},
\begin{align*}
    U_0(0,t)\mathcal{N}(\vec{u}(t))\vec{u}(t)=\begin{pmatrix}
        \tfrac{-\sin(t\jap{p})}{\jap{p}}N(u(t),x,t)u(t)\\
        \cos(t\jap{p})N(u(t),x,t)u(t)
    \end{pmatrix},
\end{align*}
using $\sin(t\jap{p})=\frac{1}{2i}\sum_{a\in\{+1,-1\}}a\,e^{iat\jap{p}}$ and $\cos(t\jap{p})=\frac{1}{2}\sum_{a\in\{+1,-1\}}e^{iat\jap{p}}$, each of the four norms in \eqref{Cor: state1} and \eqref{Cor: state2} is bounded by a constant times
\begin{equation*}
\begin{aligned}
\text{Case (i):}\quad &\sum_{a\in\{+1,-1\}}\|F_\alpha(x,p,t) e^{iat\jap{p}}\langle x\rangle^{-\sigma}\|\sup_{s\geq 0}\|\langle x\rangle^\sigma N(u(s),x,s)u(s)\|,\\
\text{Case (ii):}\quad &\sum_{a\in\{+1,-1\}}\|F_\alpha(x,p,t) e^{iat\jap{p}}\|_{L^1_x\to L^2_x}\sup_{s\geq 0}\|N(u(s),x,s)u(s)\|_{L^1_x}.
\end{aligned}
\end{equation*}
By \eqref{decay: eq1} and Assumption~\ref{asp: 2}, Case (i) lies in $L^1_t[1,\infty)$; by \eqref{decay: eq2} and Assumption~\ref{asp: 3}, so does Case (ii).
\end{proof}

\begin{lemma}\label{D: lem}
Let $\vec u(t)=(u(t),\dot u(t))$ be a global solution to the perturbed Klein-Gordon equation
\eq\label{Ppert}
\begin{cases}
(\Box+1) u(t)=N(u(t),x,t)u(t),\\
\vu(0)=(u(0),\dot{u}(0))\in \Hi,
\end{cases}\qquad (x,t)\in \mathbb{R}^n\times \mathbb{R}.
\eeq
With $\mathcal{N}(\vec u(t))$ as in \eqref{def: matrixN}, the solution $\vec u(t)$ admits the Duhamel representation
\begin{equation}\label{D: eq}
\vec u(t)=U_0(t,0)\vec u(0)+\int_0^t U_0(t,s)\,\mathcal{N}(\vec u(s))\,\vec u(s)\,ds.
\end{equation}
\end{lemma}
\begin{proof}
The system \eqref{Ppert} is equivalent to $\partial_t\vec u(t)=-(\A_0-\mathcal{N}(\vec u(t)))\vec u(t)$. Duhamel's principle applied with the free propagator $U_0(t,s)=e^{-(t-s)\A_0}$ from \eqref{Ufree} then yields \eqref{D: eq}.
\end{proof}
\subsection{Propagation estimates}\label{subsec: prop est}We adapt the propagation-observable framework introduced in \cite{SW2022} for the Schr\"odinger equation to the present matrix setting suited to the Klein-Gordon equation. The two formulations below are the vector-valued analogues, for $\vec u(t)=(u(t),\dot u(t))\in\Hi$, of the corresponding scalar estimates in \cite{SW2022}.
\begin{enumerate}
 \item(\textbf{Propagation Estimate}) Given a class of matrix operators $\{B(t)\}_{t\geq 0}$ with
 \begin{equation*}
 B(t)=\begin{pmatrix}
     b_1(t)& 0 \\
     0& b_2(t)
 \end{pmatrix},
 \end{equation*}
 we define the time-dependent inner product as:
\begin{align*}
    \langle B(t) : \vec{u}(t)\rangle_t:=(\jap{p} u_1(t),\jap{p} b_1(t) u_1(t))_{L^2_x(\mathbb{R}^n)}+( u_2(t),b_2(t)u_2(t))_{L^2_x(\mathbb{R}^n)},
\end{align*}
where $\vec{u}(t)$ denotes the solution to \eqref{KG}. Then the family $\{B(t)\}_{t\geq 0}$ is called a \textbf{Propagation Observable} if it satisfies the following condition: For a family of self-adjoint operators $B(t)$, the time derivative satisfies: there exists $L\in \mathbb{N}^+$ such that
\begin{align*}
    & \partial_t \langle B(t) : \vec{u}(t)\rangle_t=\left(\pm\sum\limits_{l=1}^L\sum\limits_{j=1}^2(\jap{p}^{2-j} u_j(t),C_{j,l}^*(t)C_{j,l}(t)\jap{p}^{2-j} u_j(t))_{L^2_x(\mathbb{R}^n)}\right)+g(t),\nonumber\\
    & g(t)\in L^1_{t}[1,\infty),\quad C_{j,l}^*(t)C_{j,l}(t)\geq0, \quad l=1,\cdots, L, j=1,2.
\end{align*}
Integrating this over time, we derive the \textbf{Propagation Estimate}: 
\begin{align*}
   \sum\limits_{l=1}^L \sum\limits_{j=1}^2\int_{t_0}^T\|C_{j,l}(t)\jap{p}^{2-j} u_j(t) \|^2dt&= \langle B(t) : \vec{u}(t)\rangle_t\vert_{t=t_0}^{t=T}-\int_{t_0}^Tg(s) ds\nonumber\\
\leq& \sup\limits_{t\in [t_0,T]} \left|\langle B(t),\vec{u}(t)\rangle_t\right|+C_g,  
\end{align*}
where $C_g:=\|g(t)\|_{L^1_t[1,\infty)}.$ 
 \item(\textbf{Relative Propagation Estimate}) Consider a class of matrix operators $\{ \tilde{B}(t)\}_{t\geq 0}$ with
 \begin{equation*}
\tilde{B}(t)=\begin{pmatrix}
     \tilde b_1(t)& 0 \\
     0& \tilde b_2(t)
 \end{pmatrix}.
 \end{equation*}
 We denote their time-dependent expectation values as:
 \begin{align*}
     \langle \tilde{B}: \vec{v}(t)\rangle_t:=(\jap{p} v_1(t), \jap{p}\tilde b_1(t) v_1(t)  )_{L^2_x(\mathbb{R}^n)}+(v_2(t), \tilde b_2(t)v_2(t)  )_{L^2_x(\mathbb{R}^n)},
 \end{align*}
 where $\vec{v}(t)$ is not necessarily the solution to \eqref{KG}, but satisfies the condition:
\eq
\sup\limits_{t\geq 0}\langle \tilde{B}:\vec{v}(t)\rangle_t<\infty. \label{phiH}
\eeq
If \eqref{phiH} holds, and if the time derivative $\partial_t\langle \tilde{B}:\vec{v}(t)\rangle_t$ meets the following estimate: there exists $L\in \mathbb{N}^+$ such that
\begin{align*}
&\partial_t\langle \tilde{B} : \vec{v}(t)\rangle_t=\left(\pm \sum\limits_{l=1}^L \sum\limits_{j=1}^2(\jap{p}^{2-j} v_{j}(t), C^*_{j,l}(t)C_{j,l}(t)\jap{p}^{2-j}v_{j}(t))_{L^2_x(\mathbb{R}^n)}\right)+g(t),\nonumber\\
&g(t)\in L^1[1,\infty), \quad C_{j,l}^*(t)C_{j,l}(t)\geq0, \quad l=1,\cdots, L, j=1,2.
\end{align*}
Then the family $\{\tilde{B}(t)\}_{t\geq 0}$ is called a {\bf Relative Propagation Observable} with respect to $\vec{v}(t)$. Integrating this over time yields the \textbf{Relative Propagation Estimate}: 
\begin{align}\label{CC}
   \sum\limits_{l=1}^L \sum\limits_{j=1}^2\int_{t_0}^T\|C_{j,l}(t)\jap{p}^{2-j} v_j(t) \|^2dt&= \langle B(t) : \vec{v}(t)\rangle_t\vert_{t=t_0}^{t=T}-\int_{t_0}^Tg(s) ds\nonumber\\
\leq& \sup\limits_{t\in [t_0,T]} \left|\langle B(t),\vec{u}(t)\rangle_t\right|+C_g.  
\end{align}
\end{enumerate}
In this paper, we set $\vec{v}(t):=U_0(0,t)\vec u(t)$ and take the operators $B_1(t),B_2(t)$ defined in \eqref{B1}, \eqref{B2}; here $F_\alpha(x,p,t)=F_c(\frac{|x|}{t^\alpha}\leq 1)F_1$ with
\eq\label{def: F1 cases}
F_1=\begin{cases}F_1(t^\beta|p|> 1)F_1(|p|\leq t^\beta) & \text{if $N$ is local (cf.\ \eqref{Falpha1}),}\\ F_1(|p|\leq t^\beta) & \text{if $N$ is non-local (cf.\ \eqref{Falpha2}),}\end{cases}
\eeq
which we abbreviate as $F_1(|p|,t)$ in the local case. We first verify that under Assumption~\ref{asp: global},
\eq
\sup\limits_{t\geq 0} \langle B_j: \vec{v}(t)\rangle_t\lesssim_{E}1,\quad j=1,2.\label{uniform: B}
\eeq
\begin{lemma}If Assumption~\ref{asp: global} holds, then \eqref{uniform: B} holds for $j=1,2$.
\end{lemma}
\begin{proof}Since the operators in $B_j(t)$ are bounded uniformly in time, $|\langle B_j:\vec{v}(t)\rangle_t|\lesssim \|\vec{v}(t)\|_\Hi^2\lesssim \|U_0(0,t)\|_{\Hi\to\Hi}^2\sup_{s\geq 0}\|\vec{u}(s)\|_\Hi^2\lesssim_E 1$.
\end{proof}
\begin{lemma}For all $\vec{v}, \vu\in \Hi$ and $j=1,2$, $(\vec{v}, B_j\vu)_{\Hi}=(B_j\vec{v}, \vu)_{\Hi}$.
\end{lemma}
\begin{proof}This follows because both $F_1F_\alpha(x,p,t)$ and $F_\alpha(x,p,t)F_c(\frac{|x|}{t^\alpha}\leq 1)$ are self-adjoint on $L^2_x(\mathbb{R}^n)$.
\end{proof}
\begin{lemma}\label{Lem: B}Given that Assumptions \ref{asp: global} and \ref{asp: 2} or \ref{asp: 3} are valid, depending on whether $\mathcal{N}$ is a local or non-local interaction respectively, both $\{B_1(t)\}_{t\geq 0}$ and $\{B_2(t)\}_{t\geq 0}$ are relative propagation observables with respect to $\vec{v}(t)$.
\end{lemma}
\begin{proof}
For simplicity, we write $F_c\equiv F_c(\frac{|x|}{t^\alpha}\leq 1)$. We compute
\begin{align*}
    \partial_t\langle B_1(t): \vec{v}(t)\rangle_t=\langle \partial_t[B_1(t)]: \vec{v}(t)\rangle_t+a_{in1}(t)+a_{in2}(t),
\end{align*}
where $a_{in1}(t),a_{in2}(t)$ collect the inner products involving $B_1(t)U_0(0,t)\mathcal{N}(\vec{u}(t))\vec{u}(t)$. By Corollary~\ref{cor: dispersive}, $a_{in1},a_{in2}\in L^1_t[1,\infty)$. Next we decompose
\begin{equation*}
\partial_t[B_1(t)]=B_{1p1}(t)+B_{1p2}(t)+B_{1r1}(t)+B_{1r2}(t),
\end{equation*}
where the $B_{1pj}$ contain $\partial_t[F_c]$ and $\partial_t[F_1]$ in symmetric (non-negative) form, namely
\begin{equation*}
B_{1p1}(t)=\begin{pmatrix}\jap{p}^{-1}F_1\partial_t[F_c]F_1 \jap{p} & 0\\ 0& F_1\partial_t[F_c]F_1\end{pmatrix}
\end{equation*}
and
\begin{equation*}
B_{1p2}(t)=2\begin{pmatrix}\jap{p} ^{-1}\sqrt{F_c}F_1\partial_t[F_1]\sqrt{F_c}\jap{p} & 0\\ 0& \sqrt{F_c}F_1\partial_t[F_1]\sqrt{F_c}\end{pmatrix},
\end{equation*}
and $B_{1r1},B_{1r2}$ are commutator remainders, each of the form $[F_1,\sqrt{F_c}]\,\partial_t[F_1]\,\sqrt{F_c}$ (or its adjoint). Since $\partial_t[F_c]\geq 0$ and $\partial_t[F_1]\geq 0$ (in Fourier space) hold in both the local and non-local cases, $\langle B_{1pj}(t):\vec{v}(t)\rangle_t\geq 0$ for $j=1,2$. For the commutator terms, Lemma~\ref{com} gives $\|[F_1,\sqrt{F_c}]\|_{L^2\to L^2}\lesssim t^{\beta-\alpha}$. Moreover, on the supports of the relevant cut-offs,
\begin{equation*}
\|\partial_t[F_1]\|\lesssim \tfrac{1}{t},
\end{equation*}
since $\partial_t F_1(t^\beta|p|)=\beta t^{\beta-1}|p|\,F_1'(t^\beta|p|)\sim \beta t^{-1}$ on the support $t^\beta|p|\sim 1$ (and analogously for $F_1(|p|\leq t^\beta)$). Combining,
\begin{equation*}
\|B_{1rj}(t)\|_{\Hi\to\Hi}\lesssim t^{\beta-\alpha-1}\in L^1_t[1,\infty)\quad\text{since } \alpha>\beta,
\end{equation*}
hence $\langle B_{1rj}(t):\vec v(t)\rangle_t\in L^1_t[1,\infty)$. Therefore $\{B_1(t)\}_{t\geq 0}$ is a relative propagation observable with respect to $\vec{v}(t)$. $\{B_2(t)\}_{t\geq 0}$ is handled similarly.
\end{proof}
\begin{corollary}\label{cor: prop}Suppose Assumption~\ref{asp: global} holds, together with either Assumption~\ref{asp: 2} (if $N$ is local) or Assumption~\ref{asp: 3} (if $N$ is non-local). Then
\eq\label{prop: eq1}
(\jap{p} v_1(t), F_1 \partial_t[F_c]F_1\jap{p} v_1(t))_{L^2_x(\mathbb{R}^n)}\in L^1_t[1,\infty),
\eeq
\eq\label{prop: eq2}
(v_2(t),  F_1\partial_t[F_c]F_1v_2(t))_{L^2_x(\mathbb{R}^n)}\in L^1_t[1,\infty),
\eeq
\eq\label{prop: eq3}
(\jap{p} v_1(t), \sqrt{F_c}F_1\partial_t[F_1] \sqrt{F_c}\jap{p} v_1(t))_{L^2_x(\mathbb{R}^n)}\in L^1_t[1,\infty)
\eeq
and 
\eq\label{prop: eq4}
(v_2(t),  \sqrt{F_c}F_1\partial_t[F_1] \sqrt{F_c}v_2(t))_{L^2_x(\mathbb{R}^n)}\in L^1_t[1,\infty),
\eeq
where $F_1$ is taken as $F_1(t^\beta|p|>1)F_1(|p|\leq t^\beta)$ when $N$ is local and as $F_1(|p|\leq t^\beta)$ when $N$ is non-local.
\end{corollary}
\begin{proof}It follows from the proof of Lemma \ref{Lem: B}.\end{proof}
\section{Proof of the existence of free channel wave operators}
We now establish, under Assumption~\ref{asp: global} together with Assumption~\ref{asp: 2} or~\ref{asp: 3}, the existence in $\Hi$ of the free channel wave operator acting on initial data $\vec u(0)$,
\begin{equation*}
\Omega_{\alpha}^*\vu(0) := s\text{-}\lim_{t\to \infty} \mathcal{F}_\alpha(x,p,t) U_0(0,t) \vu(t),
\end{equation*}
using the propagation estimates from Section~\ref{subsec: prop est}.
\begin{proof}[Proof of Theorem \ref{thm1}]Consider the vector
\begin{equation*}
\vu_{\Omega,\alpha}(t):=\mathcal{F}_\alpha(x,p,t)U_0(0,t)\vu(t).
\end{equation*}
By Cook's method, $\vu_{\Omega,\alpha}(t)=\vu_{\Omega,\alpha}(1)+\vu_{in}(t)+\vu_p(t)$, where
\begin{align*}
\vu_{in}(t)&:=-\int_1^t ds\,\mathcal{F}_\alpha(x,p,s)U_0(0,s)\mathcal{N}(\vec{u}(s))\vu(s),\\
\vu_p(t)&:=\int_1^t ds\,\partial_s[\mathcal{F}_\alpha(x,p,s)]U_0(0,s)\vu(s).
\end{align*}
Note that $\vu_{\Omega,\alpha}(1)\in \Hi$ by Assumption~\ref{asp: global} together with $\sup_{t\in\mathbb{R}}\|U_0(t,0)\|_{\Hi\to\Hi}<\infty$. Moreover, $\vu_{in}(\infty)$ exists in $\Hi$ because Corollary~\ref{cor: dispersive} gives $\|\mathcal{F}_\alpha(x,p,s)U_0(0,s)\mathcal{N}(\vec{u}(s))\vec{u}(s)\|_{\Hi}\in L^1_t[1,\infty)$.
For $\vu_p(t)$, we use the relative propagation estimates with $\vec{v}=U_0(0,t)\vec{u}(t)$ and the propagation observables $\{B_1(t)\}_{t\geq 0}, \{B_2(t)\}_{t\geq 0}$, given in \eqref{B1}, \eqref{B2}, respectively. To compute $\vu_p(t)$, we find:
\begin{align*}
    \vu_p(t)= & \int_1^t ds\mathcal{F}_c(x,p,s)\vec{v}(s)+  \int_1^t ds\mathcal{F}_1(x,p,s)\vec{v}(s)+ \int_1^t ds\mathcal{F}_r(x,p,s)\vec{v}(s)\nonumber\\
    =& \vu_{pc}(t)+\vu_{p1}(t)+\vu_{pr}(t)
\end{align*}
where $\mathcal{F}_c(x,p,s), \mathcal{F}_1(x,p,s) $ and $\mathcal{F}_r(x,p,s)$ are given by 
\begin{equation*}
\mathcal{F}_c(x,p,s)= \begin{pmatrix}
    \jap{p}^{-1} \partial_t[F_c]F_1 \jap{p} & 0\\
    0& \partial_t[F_c]F_1
\end{pmatrix},
\end{equation*}
\begin{equation*}
\mathcal{F}_1(x,p,s)= \begin{pmatrix}
    \jap{p}^{-1} \partial_t[F_1]F_c \jap{p} & 0\\
    0& \partial_t[F_1]F_c
\end{pmatrix}
\end{equation*}
and 
\begin{equation*}
\mathcal{F}_r(x,p,s)= \begin{pmatrix}
    \jap{p}^{-1} [F_c, \partial_t[F_1]] \jap{p} & 0\\
    0&[F_c,  \partial_t[F_1]]
\end{pmatrix}.
\end{equation*}
According to Lemma \ref{com}, we conclude that $\|\mathcal{F}_r(x,p,s)\vec{v}(s)\|_{\Hi}\in L^1_s[1,\infty)$. Therefore, $\vu_{pr}(\infty) $ exists in $\Hi$. \par
For $\vu_{pc}(t)=(u_{pc1}(t), u_{pc2}(t))$, applying H\"older's inequality in $s$ together with $\partial_t[F_c]\geq 0$ yields, for $t_1>t_2\geq T\geq 1$ and $j=1,2$,
\begin{align*}
\| \jap{p}^{2-j}(u_{pcj}(t_1)-u_{pcj}(t_2))\|^2 &\leq \Bigl(\int_{t_2}^{t_1}\partial_t[F_c]\,dt\Bigr)\\
&\qquad\times\int_{t_2}^{t_1}ds\,(F_1\jap{p}^{2-j}v_j(s), \partial_t[F_c]|_{t=s}F_1\jap{p}^{2-j}v_j(s))_{L^2_x}\\
&\leq \int_{t_2}^{t_1}ds\,(F_1\jap{p}^{2-j}v_j(s), \partial_t[F_c]|_{t=s}F_1\jap{p}^{2-j}v_j(s))_{L^2_x}.
\end{align*}
By \eqref{prop: eq1} and \eqref{prop: eq2} in Corollary~\ref{cor: prop}, the right-hand side tends to $0$ as $T\to\infty$, so $\{\vu_{pc}(t)\}_{t\geq 1}$ satisfies the Cauchy criterion in $\Hi$, and $\vu_{pc}(\infty)$ exists in $\Hi$. The analogous argument in the Fourier variable (using $\partial_t[\hat F_1]\geq 0$ together with \eqref{prop: eq3}, \eqref{prop: eq4}) shows that $\vu_{p1}(\infty)$ exists in $\Hi$. Combining with $\vu_{pr}(\infty)$ from above, $\vu_p(\infty)$ exists in $\Hi$. Together with $\vu_{in}(\infty)\in\Hi$, this gives $\vu_{\Omega,\alpha}(\infty)\in\Hi$ and proves the existence of $\Omega_\alpha^*\vu(0)$. The independence of $\Omega_\alpha^*\vec{u}(0)$ on $\alpha,\beta$ follows from
\begin{equation}\label{w-lim 1-Fa}
w\text{-}\lim\limits_{t\to \infty} (1-\mathcal{F}_\alpha(x,p,t))U_0(0,t)\vec{u}(t)=0\quad\text{in }\Hi,
\end{equation}
which in turn follows by testing against $\vec{v}\in C^\infty_0\times C^\infty_0$: by the definition of $F_\alpha(x,p,t)$ (see \eqref{Falpha1}, \eqref{Falpha2}), each cutoff factor tends pointwise to $1$ as $t\to\infty$, hence $w\text{-}\!\lim_{t\to\infty}(1-F_\alpha(x,p,t))=0$ on $L^2_x(\mathbb{R}^n)$; consequently $\|(1-\mathcal{F}_\alpha(x,p,t))\vec v\|_\Hi\to 0$ for every such $\vec v$.
\end{proof}
\begin{proof}[Proof of Theorem \ref{thm2}]We follow the same argument as in the proof of Theorem~\ref{thm1}, with the local choice $F_1=F_1(t^\beta|p|> 1)F_1(|p|\leq t^\beta)$ in \eqref{def: F1 cases} replaced by the non-local choice $F_1=F_1(|p|\leq t^\beta)$, and with case~(ii) of Corollary~\ref{cor: dispersive} used in place of case~(i). The propagation observables $\{B_j(t)\}_{t\geq 0}$, $j=1,2$, defined in \eqref{B1}, \eqref{B2}, remain relative propagation observables with respect to $\vec{v}(t)=U_0(0,t)\vec{u}(t)$, by Lemma~\ref{Lem: B}. The independence of $\Omega_\alpha^*\vec u(0)$ on $(\alpha,\beta)$ follows from the same weak-limit argument as in the proof of Theorem~\ref{thm1}.
\end{proof}

\section{Proof of properties of the weakly localized part}
In this section we prove Theorem~\ref{thm3}, that is, the sub-ballistic estimate \eqref{conlusion: wlocal} for the weakly localized part $\vu_{wlc}(t)$. The proof rests on a phase-space decomposition into \emph{forward} and \emph{backward propagation waves}, which we introduce next.
\subsection{Forward/backward propagation waves}\label{subsec: fbw}
We follow the forward/backward propagation framework of \cite[Section~4.1]{SW2022} (where these waves are introduced as the spatial analogue of Enss's incoming/outgoing waves \cite{M1979}); the only adjustment is in the final choice of group velocity, specialized here to the Klein-Gordon case $v(p)=p/\jap{p}$. We recall the relevant ingredients below.

Fix a smooth partition of unity $\{F^{\hat h}\}_{\hat h\in I}$ on the unit sphere $S^{n-1}$, with index set
\begin{equation*}
I=\{\hat h_1,\ldots,\hat h_N\}\subseteq S^{n-1},
\end{equation*}
and a companion smooth cutoff $\tilde F^{\hat h}:S^{n-1}\to\mathbb{R}$, both built from a parameter $c\in (0,\tfrac{1}{200})$ via
\begin{equation*}
F^{\hat h}(\xi)=\begin{cases}1 & |\xi-\hat h|<c\\ 0 & |\xi-\hat h|>2c\end{cases},\qquad
\tilde F^{\hat h}(\xi)=\begin{cases}1 & |\xi-\hat h|<4c\\ 0 & |\xi-\hat h|>8c\end{cases},\qquad \xi\in S^{n-1},
\end{equation*}
chosen so that (with $\hat h:=h/|h|$ for $h\in\mathbb{R}^n\setminus\{0\}$) the geometric inequalities
\eq\label{Feq1}
F^{\hat h}(\hat x)\tilde F^{\hat h}(\hat q)\,|x+q|\geq \tfrac{1}{10}(|x|+|q|),
\eeq
\eq\label{Feq2}
F^{\hat h}(\hat x)(1-\tilde F^{\hat h}(\hat q))\,|x-q|\geq \tfrac{1}{10^6}(|x|+|q|)
\eeq
hold for all $x,q\in\mathbb{R}^n$, with the convention $\hat 0 := 0$. The corresponding projections on the forward/backward propagation set (with respect to $(r,v)\in\mathbb{R}^{n+n}$) are
\eq\label{Prv+}
P^\pm(r,v):=\sum_{b=1}^N F^{\hat h_b}(\hat r)\,\tilde F^{\hat h_b}_\pm(\hat v),\qquad \tilde F^{\hat h_b}_+:=\tilde F^{\hat h_b},\quad \tilde F^{\hat h_b}_-:=1-\tilde F^{\hat h_b}.
\eeq
We further write $F_a(\lambda):=F(\lambda>a):=F(\lambda/a)$ for the rescaled smooth cutoff with $F(k)=1$ when $k\geq 1$ and $F(k)=0$ when $k\leq \tfrac{1}{2}$.

For the Klein-Gordon flow we specialize to
\begin{equation*}
P^\pm\equiv P^\pm(|x|,v(p)),\qquad v(p):=\nabla_p\omega(p)=\tfrac{p}{\sqrt{|p|^2+1}},\quad \omega(p):=\sqrt{|p|^2+1},
\end{equation*}
so that $v(p)$ is the KG group velocity. This is the only place where the dispersion relation enters; all subsequent geometric estimates rely on \eqref{Feq1}, \eqref{Feq2}, \eqref{Prv+}.
\begin{lemma}[Non-zero-velocity propagation estimate]\label{Lem: Pprop}For all $\sigma>0$, $1>\alpha>\beta>0$, $s\geq 0$ and $t\geq 1$, the operator norm bound
\eq\label{Lemmain:goal}
\|F(|x|> t^\alpha)\,P^\pm e^{\pm is\omega(p)}\,F_{t^{-\beta}}(v(p))\,\langle x\rangle^{-\sigma}\|\lesssim_{n,\sigma} \frac{1}{(t^\alpha+st^{-\beta})^\sigma}
\eeq
holds, where $\omega(p)=\sqrt{|p|^2+1}$.
\end{lemma}
\begin{proof}Fix $f\in L^2_x(\mathbb{R}^n)$. We assume that
\eq\label{f cutoff}
f=\chi(|x|<\tfrac{1}{10^8}(t^\alpha+st^{-\beta}))\,f,
\eeq
since $\langle x\rangle^{-\sigma}\lesssim_\sigma 1/(t^\alpha+st^{-\beta})^\sigma$ on $\{|x|\geq\tfrac{1}{10^8}(t^\alpha+st^{-\beta})\}$, so the complementary contribution
\begin{equation*}
\bigl\|F(|x|>t^\alpha)P^\pm e^{\pm is\omega(p)}F_{t^{-\beta}}(v(p))\langle x\rangle^{-\sigma}\bigl(1-\chi(|x|<\tfrac{1}{10^8}(t^\alpha+st^{-\beta}))\bigr)f\bigr\|\lesssim_\sigma\tfrac{\|f\|}{(t^\alpha+st^{-\beta})^\sigma}
\end{equation*}
is already of the required order. Set
\begin{equation*}
Q(t,s)f:=F(|x|> t^\alpha)\,P^\pm e^{\pm is\omega(p)}\,F_{t^{-\beta}}(v(p))\,\langle x\rangle^{-\sigma}f.
\end{equation*}
Let $\{F_\xi\}_{\xi\in\mathbb{Z}^n}$ be a smooth Littlewood--Paley partition of unity on $\mathbb{R}^n$ with
\begin{equation*}
F_\xi(x)=\begin{cases}1 & |x-\xi|\leq 1,\\ 0 & |x-\xi|>2,\end{cases}
\end{equation*}
and decompose
\begin{align*}
Q(t,s)f&=\sum_{\eta=(\eta_1,\eta_2,\eta_3)\in\mathbb{Z}^{3n}}Q_\eta(t,s)f,\\
Q_\eta(t,s)f&:=F_{\eta_1}(x)F(|x|>t^\alpha)P^\pm e^{\pm is\omega(p)}F_{\eta_2}(v(p))F_{t^{-\beta}}(v(p))\langle x\rangle^{-\sigma}F_{\eta_3}(x)f.
\end{align*}
Using the Fourier inversion formula \eqref{eq: Fourier} together with \eqref{Prv+},
\begin{align*}
     Q_\eta(t,s)f=&F_{\eta_1} F(|x|> t^\alpha) P^\pm e^{\pm i s\omega(p)} F_{\eta_2}(v(p))F_{t^{-\beta}}(v(p))\langle x\rangle^{-\sigma}F_{\eta_3}(x)f\nonumber\\
     =&  c_n \int_{\R^{2n}} F_{\eta_1}(x)F(|x|> t^\alpha) P^\pm e^{\pm i s\omega(q)} e^{-iq\cdot y}F_{\eta_2}(q)F_{t^{-\beta}}(v(q))\nonumber\\
     &\qquad\qquad\qquad\qquad\times F_{\eta_3}(y)\langle y\rangle^{-\sigma} f(y)d^nq\,d^ny\nonumber\\
     =& \sum\limits_{b=1}^N  c_n \int_{\R^{2n}} F_{\eta_1}(x)F(|x|> t^\alpha)  F^{\hat{h}_b}(\hat{x})\tilde{F}^{\hat{h}_b}_\pm(\widehat{v(q)})e^{\pm i s\omega(q)} e^{-iq\cdot y}F_{\eta_2}(q)F_{t^{-\beta}}(v(q))\nonumber\\
     &\qquad\qquad\qquad\qquad\times F_{\eta_3}(y)\langle y\rangle^{-\sigma} f(y)d^nq\,d^ny
\end{align*}
On the support of $P^\pm(r,v(q))F_{t^{-\beta}}(v(q))$, the geometric inequalities \eqref{Feq1}--\eqref{Feq2}, the spatial cutoff $|x|>t^\alpha$, the velocity bound $|\eta_2|>t^{-\beta}/2$ (from the support of $F_{\eta_2}(v(p))F_{t^{-\beta}}(v(p))$), and \eqref{f cutoff} (which gives $|\eta_3|<\tfrac{1}{10^8}(t^\alpha+st^{-\beta})$) together yield, for all $s\geq 0$,
\begin{equation}\label{geom: Pprop}
|\eta_1\pm s\eta_2-\eta_3|\geq \tfrac{1}{10^6}(|\eta_1|+s|\eta_2|)-|\eta_3|-\tfrac{1}{10^7}(t^\alpha+st^{-\beta})\geq \tfrac{1}{10^7}(|\eta_1|+|\eta_3|+t^\alpha+st^{-\beta}).
\end{equation}
We define $\tilde{F}_\xi$ as a smooth cut-off function satisfying 
\begin{equation*}
    \tilde{F}_\xi(x)=\begin{cases}
        1 &\text{ if }|x-\xi|\leq 10\\
        0 & \text{ if }|x-\xi|>20
    \end{cases}.
\end{equation*}
Applying the method of non-stationary phase ($M+2n+2$ integrations by parts in $q$), together with Plancherel's theorem and the estimates
\begin{equation*}
\| F_{\eta_1}(x)F^{\hat{h}_b}(\hat{x})g\|\leq \|g\|\quad\qquad \forall g\in L^2_x(\mathbb{R}^n),
\end{equation*}
\begin{equation*}
\|\tilde{F}_{\xi}(v(p)) F_{\eta_1}(x)^2F(|x|> t^\alpha)\tilde{F}_{\xi'}(v(p)) \|\lesssim_n \frac{1}{\langle \xi-\xi'\rangle^{n+1}},
\end{equation*}
\begin{equation*}
    \| \tilde{F}^{\hat{h}_b}_\pm(\widehat{v(q)}) F_{t^{-\beta}}(v(q))F_{\eta_2}(q)\|_{C^{M+n+1}_q(\mathbb{R}^n)}\lesssim_{M+n} t^{(M+n+1)\beta}
\end{equation*}
we obtain
\begin{align*}
\|\sum_{\eta_2\in \mathbb{Z}^n}Q_\eta(t,s)f\|^2\lesssim_{n,M} & \tfrac{1}{\langle |\eta_1|+|\eta_3|\rangle^{4(n+1)}}\tfrac{1}{(t^\alpha+st^{-\beta})^{2M}}t^{(M+n+1)\beta}\nonumber\\
&\times\sum_{\eta_2,\eta_2'\in\mathbb{Z}^n}\tfrac{\|\tilde F_{\eta_2}(v(p))F_{\eta_3}(x)\langle x\rangle^{-\sigma}f\|\,\|\tilde F_{\eta_2'}(v(p))F_{\eta_3}(x)\langle x\rangle^{-\sigma}f\|}{\langle \eta_2-\eta_2'\rangle^{n+1}}\nonumber\\
\lesssim_{n,M} & \tfrac{1}{\langle |\eta_1|+|\eta_3|\rangle^{4(n+1)}}\tfrac{1}{(t^\alpha+st^{-\beta})^{2M}}t^{(M+n+1)\beta}\|F_{\eta_3}(x)\langle x\rangle^{-\sigma}f\|^2,
\end{align*}
where in the last step we have used
\begin{equation*}
\sum\limits_{\eta_2\in\mathbb{Z}^n}\sum\limits_{\eta_2'\in\mathbb{Z}^n}\tfrac{1}{\langle \eta_2-\eta_2'\rangle^{n+1}}\|\tilde F_{\eta_2}(v(p))F_{\eta_3}(x)\langle x\rangle^{-\sigma}f\|\,\|\tilde F_{\eta_2'}(v(p))F_{\eta_3}(x)\langle x\rangle^{-\sigma}f\|\lesssim \|F_{\eta_3}(x)\langle x\rangle^{-\sigma}f\|^2.
\end{equation*}
Taking square roots and summing over $(\eta_1,\eta_3)\in\mathbb{Z}^{2n}$, then choosing $M$ large enough that $\alpha M-(M+n+1)\beta\geq\alpha\sigma$ (recall $\alpha>\beta$), we conclude
\begin{equation*}
\|Q(t,s)f\|\leq \sum_{(\eta_1,\eta_3)\in\mathbb{Z}^{2n}}\|\sum_{\eta_2\in\mathbb{Z}^n}Q_\eta(t,s)f\|\lesssim_{n,\sigma}\tfrac{1}{(t^\alpha+st^{-\beta})^\sigma}\|\langle x\rangle^{-\sigma}f\|.\qedhere
\end{equation*}
\end{proof}
\begin{lemma}[Backward-projection propagation estimate]\label{Lem: Pprop2}For all $\sigma>0$, $\beta\in(0,1)$ and $t\geq 1$, the operator norm bound
\begin{equation*}
\|P^\pm e^{\pm it\omega(p)}\,F_{t^{-\beta}}(v(p))\,\langle x\rangle^{-\sigma}\|\lesssim_{n,\sigma}\frac{1}{(t^{1-\beta})^\sigma}
\end{equation*}
holds.
\end{lemma}
\begin{proof}The proof is completely analogous to that of Lemma~\ref{Lem: Pprop} (with $F(|x|>t^\alpha)$ dropped, $s$ replaced by $t$, and $t^\alpha$ removed from the right-hand side; in particular the cutoff \eqref{f cutoff} becomes $f=\chi(|x|<\tfrac{1}{10^8}t^{1-\beta})f$). The geometric inequality \eqref{geom: Pprop} is replaced by
\begin{equation*}
|\eta_1\pm t\eta_2-\eta_3|\geq \tfrac{1}{10^6}(|\eta_1|+t|\eta_2|)-|\eta_3|-\tfrac{1}{10^7}t^{1-\beta}\geq \tfrac{1}{10^7}(|\eta_1|+|\eta_3|+t^{1-\beta}),
\end{equation*}
which follows from \eqref{Feq1}, \eqref{Feq2} and $|\eta_3|<\tfrac{1}{10^8}t^{1-\beta}$ alone, without using $|x|>t^\alpha$.
\end{proof}
\begin{corollary}[Zero-velocity propagation estimate]\label{cor: small}Let $\sigma$ be as in Lemma~\ref{Lem: Pprop}. If $s\in [0,t]$, $\alpha>\frac{1}{2}$ and $\beta=\frac{1}{2}$, then
\eq\label{cor: eq}
\|F(|x|> t^\alpha)(1-F_{t^{-\beta}}(v(p)))e^{\pm is\omega(p)}\langle x\rangle^{-\sigma}\|\leq \frac{C}{t^{\alpha\sigma}}
\eeq
for some constant $C=C(n,\sigma)>0$.
\end{corollary}
\begin{proof}Fix $f\in L^2_x(\mathbb{R}^n)$ and set
\begin{equation*}
    \tilde Q(t,s)f:=F(|x|> t^\alpha)e^{\pm is\omega(p)}(1-F_{t^{-\beta}}(v(p)))\langle x\rangle^{-\sigma}f.
\end{equation*}
Since
\begin{align*}
\|\langle x\rangle^{-\sigma}\chi(|x|>\tfrac{1}{10^7}(t^\alpha-st^{-\beta}))\|\lesssim \tfrac{1}{t^{\alpha\sigma}},
\end{align*}
it suffices to consider $f\in L^2_x(\mathbb{R}^n)$ supported in $\{|x|\leq\tfrac{1}{10^7}(t^\alpha-st^{-\beta})\}$. Decomposing via Littlewood--Paley with $\bar F_{t^{-\beta}}:=1-F_{t^{-\beta}}$ and $F_\xi$ the unit-scale partition,
\begin{align*}
\tilde Q(t,s)f&=\sum_{\eta=(\eta_1,\eta_2,\eta_3)\in\mathbb{Z}^{3n}}\tilde Q_\eta(t,s)f,\\
\tilde Q_\eta(t,s)f&:=F_{\eta_1}(x)F(|x|> t^\alpha)e^{\pm is\omega(p)}F_{\eta_2}(v(p))\bar F_{t^{-\beta}}(v(p))\langle x\rangle^{-\sigma}F_{\eta_3}(x)f.
\end{align*}
The support assumption on $f$ together with $|\eta_1|\geq t^\alpha/2$ gives $|\eta_3|\leq 2|\eta_1|/10^7$; combined with $\alpha+\beta>1$ and $s\in [0,t]$ (so $t^\alpha\geq t^{1-\beta}\geq st^{-\beta}$), this yields, for $j\in\{+1,-1\}$,
\begin{equation*}
|\eta_1+js\eta_2-\eta_3|\geq \tfrac{1}{10^6}(|\eta_1|-st^{-\beta})-\tfrac{1}{10^7}(t^\alpha+st^{-\beta})\geq \tfrac{1}{10^7}(|\eta_1|+|\eta_3|+t^\alpha).
\end{equation*}
Moreover, each $p$-derivative of $\bar F_{t^{-\beta}}(v(p))$ produces a factor of order $t^\beta$, which combined with the oscillation factor $t^\alpha$ from the geometric inequality above yields a net decay $t^{\beta-\alpha}$ per integration by parts (since $\alpha>\beta$).
Applying the method of non-stationary phase via the same Littlewood--Paley and Plancherel argument as in the proof of Lemma~\ref{Lem: Pprop} (the only change being the substitution of $\bar F_{t^{-\beta}}(v(p))$ for $F_{t^{-\beta}}(v(p))$ throughout) then yields \eqref{cor: eq}.
\end{proof}
\begin{remark}\label{rem: cor 4.3 sharp}The condition $\alpha>\tfrac{1}{2}$ in Corollary~\ref{cor: small} is essentially sharp, and is the reason the same lower bound on the spatial-cutoff exponent appears in Theorem~\ref{thm3}. Indeed, for $\alpha\leq\tfrac{1}{2}$ no choice of $\beta$ closes the argument:
\begin{itemize}
\item If $\beta\geq\tfrac{1}{2}$, the non-stationary phase method breaks down: each $p$-derivative of $F_{t^{-\beta}}(v(p))$ produces a factor of order $t^\beta$, to be compensated by the oscillation factor $t^\alpha$; for the compensation to leave any decay, one needs $\alpha>\beta$, which fails when $\beta\geq\tfrac{1}{2}\geq\alpha$.
\item If $\beta<\tfrac{1}{2}$, the geometric inequality used in the support reduction, namely $|\eta_3|\leq\tfrac{1}{10^7}(t^\alpha-st^{-\beta})$, fails because $st^{-\beta}$ may exceed $t^\alpha$ for $s\in[0,t]$.
\end{itemize}
Hence $\alpha>\tfrac{1}{2}$ is forced; the choice $\beta=\tfrac{1}{2}$ is a feasible balance that makes both estimates simultaneously available.
\end{remark}

\begin{lemma}\label{lem: low freq conv}For every $f\in L^2_x(\mathbb{R}^n)$,
\begin{equation*}
\|(1-F_{t^{-\beta}}(v(p)))f\|\to 0\quad\text{as }t\to\infty.
\end{equation*}
\end{lemma}
\begin{proof}Since $F_{t^{-\beta}}(v(p))=F(|v(p)|>t^{-\beta})\to 1$ a.e.\ in $p$-space as $t\to\infty$ (because $|v(p)|>0$ a.e.) and $|(1-F_{t^{-\beta}}(v(p)))f|\leq|f|$, the result follows by the dominated convergence theorem.
\end{proof}
\subsection{Proof of Theorem~\ref{thm3}}We outline the proof of Theorem~\ref{thm3}. A caveat on notation: throughout this subsection the parameters $\alpha,\beta$ are re-used to denote new quantities subject to $\alpha>\max\{\tfrac{1}{2},\tfrac{1}{\sigma}\}$ and $\beta=\tfrac{1}{2}$ (as required by Lemma~\ref{Lem: Pprop} and Corollary~\ref{cor: small}; see Remark~\ref{rem: cor 4.3 sharp}). These constraints are disjoint from the admissible range $\mathscr{I}_{\mathrm{loc}}$ of Theorem~\ref{thm1}. This is harmless because, by Theorem~\ref{thm1}, $\Omega_\alpha^*\vec u(0)\in\Hi$ exists and is independent of $(\alpha,\beta)\in\mathscr{I}_{\mathrm{loc}}$; we invoke Theorem~\ref{thm1} once with any admissible pair to fix the asymptotic state and then work with the new $\alpha,\beta$ below. Write $U_0(t,0)\Omega_\alpha^*\vec u(0)$ as
\begin{equation*}
U_0(t,0)\Omega_\alpha^*\vec{u}(0)=: (u_{\Omega,1}(t),u_{\Omega,2}(t)).
\end{equation*}
In order to prove \eqref{conlusion: wlocal}, it suffices to show that for all $\alpha>\max\{\tfrac{1}{2},\tfrac{1}{\sigma}\}$, as $t\to \infty$, both
\eq\label{weak: goal1}
\| F(|x|>t^\alpha)\left(\jap{p} u_1(t)- \jap{p} u_{\Omega,1}(t)\right)\|\to 0
\eeq
and
\eq\label{weak: goal2}
\| F(|x|>t^\alpha)\left( u_2(t)- u_{\Omega,2}(t)\right)\|\to 0
\eeq
hold true, then by taking
\begin{equation*}
\vec{u}_{wlc}(t):=F(|x|\leq t^\alpha)\vec u(t),
\end{equation*}
\eqref{weak: goal1} and \eqref{weak: goal2}, together with \eqref{D: eq} and the fact that $\|F(|x|\leq t^\alpha)U_0(t,0)\Omega_\alpha^*\vec u(0)\|_\Hi\to 0$ as $t\to\infty$ (the free wave concentrates on $|x|\sim t\gg t^\alpha$), imply \eqref{limit} and \eqref{conlusion: wlocal}. It therefore remains to prove \eqref{weak: goal1} and \eqref{weak: goal2}. To this end, by \eqref{D: eq} in Lemma~\ref{D: lem}, $u(t)$ and $\dot u(t)$ read
\begin{align}\label{ut: form}
    u(t)=& \cos(t\jap{p})u(0)+\frac{\sin(t\jap{p})}{\jap{p}}\dot u(0) \nonumber\\
    &+\int_0^t  \frac{\sin((t-s)\jap{p})}{ \jap{p}}N(u(s),x,s)u(s) ds
\end{align}
and
\begin{align}\label{dotut: form}
    \dot u(t)=&-\sin(t\jap{p})\jap{p}u(0)+\cos(t\jap{p})\dot u(0)\nonumber\\
    &+\int_0^t \cos((t-s)\jap{p})N(u(s),x,s)u(s).
\end{align}
Based on \eqref{ut: form} and \eqref{dotut: form}, we write $u(t)$ and $\dot u(t)$ {in terms of $e^{\pm it\jap{p}}$ flows:}
\begin{align*}
     u(t)=&\frac{1}{2}\left(e^{it\jap{p}}+e^{-it\jap{p}}\right)u(0)+\frac{1}{2i\jap{p}}\left( e^{it\jap{p}}-e^{-it\jap{p}}\right)\dot u(0)\nonumber\\
     &+\int_0^t \frac{1}{2i\jap{p}} \left( e^{i(t-s)\jap{p}}-e^{-i(t-s)\jap{p}}\right)N(u(s),x,s)u(s)ds
\end{align*}
and
\begin{align*}
     \dot u(t)=&\frac{-1}{2i}\left(e^{it\jap{p}}-e^{-it\jap{p}}\right)\jap{p}u(0)+\frac{1}{2}\left( e^{it\jap{p}}+e^{-it\jap{p}}\right)\dot u(0)\nonumber\\
     &+\int_0^t \frac{1}{2} \left( e^{i(t-s)\jap{p}}+e^{-i(t-s)\jap{p}}\right)N(u(s),x,s)u(s)ds.
\end{align*}
That is,
\eq\label{decomp: u total}
\begin{cases}
u(t)=&u(t,+)+u(t,-)\\
\dot u(t)=&\dot u(t,+)+\dot u(t,-)
\end{cases},
\eeq
where $u(t,+), u(t,-), \dot u(t,+)$ and $\dot u(t,-)$ are defined as:
\begin{align}\label{ut+}
    u(t,+):=& \frac{1}{2}e^{it\jap{p}}u(0)+\frac{1}{2i\jap{p}} e^{it\jap{p}}\dot u(0)\nonumber\\
     &+\int_0^t \frac{1}{2i\jap{p}} e^{i(t-s)\jap{p}}N(u(s),x,s)u(s)ds,
\end{align}
\begin{align}\label{ut-}
    u(t,-):=& \frac{1}{2}e^{-it\jap{p}}u(0)-\frac{1}{2i\jap{p}}e^{-it\jap{p}}\dot u(0)\nonumber\\
     &-\int_0^t \frac{1}{2i\jap{p}} e^{-i(t-s)\jap{p}}N(u(s),x,s)u(s)ds,
\end{align}
\begin{align*}
   \dot u(t,+):=&\frac{-1}{2i}e^{it\jap{p}}\jap{p}u(0)+\frac{1}{2} e^{it\jap{p}}\dot u(0)\nonumber\\
     &+\int_0^t \frac{1}{2}  e^{i(t-s)\jap{p}}N(u(s),x,s)u(s)ds
\end{align*}
and
\begin{align}\label{dotut-}
    \dot u(t,-):=& \frac{1}{2i}e^{-it\jap{p}}\jap{p}u(0)+\frac{1}{2}e^{-it\jap{p}}\dot u(0)\nonumber\\
     &+\int_0^t \frac{1}{2} e^{-i(t-s)\jap{p}}N(u(s),x,s)u(s)ds.
\end{align}
To prove Theorem~\ref{thm3}, we require the following lemmas.
\begin{lemma}\label{lem: omega exist}Assume all assumptions in Theorem~\ref{thm1} are satisfied. Then $\jap{p}\,u(t,\pm)$ and $\dot u(t,\pm)$ are bounded in $L^2_x(\mathbb{R}^n)$ uniformly in time, and the weak limits
\eq\label{uo pm}
\jap{p}\,u_{\Omega,\pm}:=w\text{-}\!\lim_{t\to\infty}e^{\mp it\jap{p}}\jap{p}\,u(t,\pm),
\eeq
\eq\label{duo pm}
\dot u_{\Omega,\pm}:=w\text{-}\!\lim_{t\to\infty}e^{\mp it\jap{p}}\dot u(t,\pm)
\eeq
exist in $L^2_x(\mathbb{R}^n)$.
\end{lemma}
\begin{proof}A direct computation from \eqref{ut+}--\eqref{dotut-} yields the key identity
\eq\label{key id}
\dot u(t,\pm)=\pm i\jap{p}\,u(t,\pm),
\eeq
which, combined with $u(t)=u(t,+)+u(t,-)$ and $\dot u(t)=\dot u(t,+)+\dot u(t,-)$, gives
\eq\label{ut+-: solve}
\jap{p}\,u(t,\pm)=\tfrac{1}{2}\jap{p}\,u(t)\mp\tfrac{i}{2}\dot u(t),\qquad \dot u(t,\pm)=\pm\tfrac{i}{2}\jap{p}\,u(t)+\tfrac{1}{2}\dot u(t).
\eeq
Since $\|\vec u(t)\|_\Hi$ is bounded uniformly in time by Assumption~\ref{asp: global}, \eqref{ut+-: solve} yields the same uniform-in-time bound for $\jap{p}\,u(t,\pm)$ and $\dot u(t,\pm)$ in $L^2_x(\mathbb{R}^n)$.

By Theorem~\ref{thm1} and \eqref{w-lim 1-Fa}, $w\text{-}\!\lim_{t\to\infty}U_0(0,t)\vec u(t)=\Omega_\alpha^*\vec u(0)=:(u_{\Omega,1},u_{\Omega,2})$ exists in $\Hi$; in particular, the first component admits a weak limit in $H^1(\mathbb{R}^n)$ and the second in $L^2_x(\mathbb{R}^n)$. A direct computation using \eqref{key id} and $\cos(t\jap{p})\mp i\sin(t\jap{p})=e^{\mp it\jap{p}}$ shows
\eq\label{U0: split}
U_0(0,t)\vec u(t)=\bigl(e^{-it\jap{p}}u(t,+)+e^{it\jap{p}}u(t,-),\;e^{-it\jap{p}}\dot u(t,+)+e^{it\jap{p}}\dot u(t,-)\bigr).
\eeq
Applying $\jap{p}$ to the first component of \eqref{U0: split} (so that its weak limit $\jap{p}\,u_{\Omega,1}$ lies in $L^2_x(\mathbb{R}^n)$), combining with the weak-limit version of \eqref{key id}, $\dot u_{\Omega,\pm}=\pm i\jap{p}\,u_{\Omega,\pm}$, and solving the resulting $2\times 2$ linear systems for $\jap{p}\,u_{\Omega,\pm}$ and $\dot u_{\Omega,\pm}$ yields the existence of the weak limits \eqref{uo pm}--\eqref{duo pm} in $L^2_x(\mathbb{R}^n)$, together with the linear-combination representation
\eq\label{uOmega: lc}
\jap{p}\,u_{\Omega,\pm}=\tfrac{1}{2}\jap{p}\,u_{\Omega,1}\mp\tfrac{i}{2}u_{\Omega,2},\qquad \dot u_{\Omega,\pm}=\pm\tfrac{i}{2}\jap{p}\,u_{\Omega,1}+\tfrac{1}{2}u_{\Omega,2}.
\eeq
Their boundedness in $L^2_x(\mathbb{R}^n)$ then follows from \eqref{uOmega: lc}, since $\jap{p}\,u_{\Omega,1}$ and $u_{\Omega,2}$ both lie in $L^2_x(\mathbb{R}^n)$ by $\Omega_\alpha^*\vec u(0)\in\Hi$.
\end{proof}
Recall that the operator $P^\pm$ are given by 
\begin{equation*}
P^\pm =P^\pm(r,v(p)),\quad r=|x|,\quad v(p)=\frac{p}{\sqrt{1+|p|^2}}.
\end{equation*}
We define the projections on the forward/backward propagation set as
\begin{equation*}
P_+^\pm\equiv P^{\mp},\qquad P_-^\pm\equiv P^{\pm}.
\end{equation*}
Here the superscript $+$ (resp.~$-$) corresponds to the forward (resp.~backward) propagation set, while the subscript $\pm$ specifies the flow $e^{it(\pm\jap{p})}$; that is, $P^+_\pm$ and $P^-_\pm$ are the projections onto the forward and backward propagation sets, respectively, for the flow $e^{it(\pm\jap{p})}$. We decompose $\jap{p}u(t,\pm)$ and $\dot u(t,\pm)$ first using the velocity cutoff $F_{t^{-\beta}}(v(p))$ into non-zero- and zero-velocity parts, and then split the non-zero-velocity part into forward/backward waves via $P^\pm_\pm$:
\begin{align}\label{decomp: u sub}
\jap{p}u(t,\pm)=& P^+_\pm F_{t^{-\beta}}(v(p))\jap{p}u(t,\pm)+P^-_\pm F_{t^{-\beta}}(v(p))\jap{p}u(t,\pm)\nonumber\\
&+(1-F_{t^{-\beta}}(v(p)))\jap{p}u(t,\pm),
\end{align}
\begin{align*}
\dot u(t,\pm)=& P^+_\pm F_{t^{-\beta}}(v(p))\dot u(t,\pm)+P^-_\pm F_{t^{-\beta}}(v(p))\dot u(t,\pm)\nonumber\\
&+(1-F_{t^{-\beta}}(v(p)))\dot u(t,\pm).
\end{align*}
The forward, non-zero-velocity pieces are approximated by the corresponding propagated $u_{\Omega,\pm}$, $\dot u_{\Omega,\pm}$:
\begin{align*}
P^+_\pm F_{t^{-\beta}}(v(p))\jap{p}u(t,\pm)=&P^+_\pm F_{t^{-\beta}}(v(p))e^{\pm it\jap{p}}\jap{p}u_{\Omega,\pm}\nonumber\\
&+\bigl(P^+_\pm F_{t^{-\beta}}(v(p))\jap{p}u(t,\pm)-P^+_\pm F_{t^{-\beta}}(v(p))e^{\pm it\jap{p}}\jap{p}u_{\Omega,\pm}\bigr),
\end{align*}
and similarly for $\dot u(t,\pm)$ in terms of $\dot u_{\Omega,\pm}$.

The backward, non-zero-velocity pieces are approximated by the free parts of \eqref{ut+}--\eqref{dotut-}, namely by $e^{\pm it\jap{p}}\jap{p}u(0,\pm)$ and $e^{\pm it\jap{p}}\dot u(0,\pm)$, where
\begin{equation*}
u(0,\pm):=\tfrac{1}{2}u(0)\pm\tfrac{1}{2i\jap{p}}\dot u(0),\qquad \dot u(0,\pm):=\mp\tfrac{1}{2i}\jap{p}u(0)+\tfrac{1}{2}\dot u(0)
\end{equation*}
denote the initial values of $u(t,\pm)$ and $\dot u(t,\pm)$:
\begin{align*}
P^-_\pm F_{t^{-\beta}}(v(p))\jap{p}u(t,\pm)=&P^-_\pm F_{t^{-\beta}}(v(p))e^{\pm it\jap{p}}\jap{p}u(0,\pm)\nonumber\\
&+\bigl(P^-_\pm F_{t^{-\beta}}(v(p))\jap{p}u(t,\pm)-P^-_\pm F_{t^{-\beta}}(v(p))e^{\pm it\jap{p}}\jap{p}u(0,\pm)\bigr),
\end{align*}
and similarly for $\dot u(t,\pm)$ in terms of $\dot u(0,\pm)$.

In Lemma~\ref{lem: weak err} below, we apply Lemma~\ref{Lem: Pprop} to the non-zero-velocity pieces and Corollary~\ref{cor: small} to the zero-velocity piece to prove the following error estimates: as $t\to\infty$, for every $\alpha>\max\{\tfrac{1}{2},\tfrac{1}{\sigma}\}$,
\eq\label{weak: eq1}
\bigl\| F(|x|>t^\alpha)\,P^+_\pm F_{t^{-\beta}}(v(p))\bigl(\jap{p}u(t,\pm)-e^{\pm it\jap{p}}\jap{p}u_{\Omega,\pm}\bigr)\bigr\|\to 0,
\eeq
\eq\label{weak: eq2}
\bigl\| F(|x|>t^\alpha)\,P^+_\pm F_{t^{-\beta}}(v(p))\bigl(\dot u(t,\pm)-e^{\pm it\jap{p}}\dot u_{\Omega,\pm}\bigr)\bigr\|\to 0,
\eeq
\eq\label{weak: eq3}
\bigl\| F(|x|>t^\alpha)\,P^-_\pm F_{t^{-\beta}}(v(p))\,\xi(t)\bigr\|\to 0,
\eeq
and
\eq\label{weak: eq4}
\bigl\| F(|x|>t^\alpha)\,(1-F_{t^{-\beta}}(v(p)))\,\xi(t)\bigr\|\to 0
\eeq
for every $\xi(t)\in\bigl\{\jap{p}u(t,\pm),\,\dot u(t,\pm)\bigr\}$.
\begin{lemma}\label{lem: weak err}If Assumption~\ref{asp: 2} is satisfied for some $\sigma>2$, then \eqref{weak: eq1}--\eqref{weak: eq4} hold for every $\alpha>\max\{\tfrac{1}{2},\tfrac{1}{\sigma}\}$.
\end{lemma}
\begin{proof}We prove \eqref{weak: eq1}, \eqref{weak: eq3} and \eqref{weak: eq4} for $\xi(t)=\jap{p}u(t,\pm)$; \eqref{weak: eq2}, and \eqref{weak: eq3}--\eqref{weak: eq4} for $\xi(t)=\dot u(t,\pm)$, follow analogously.

\noindent\fbox{\emph{Proof of \eqref{weak: eq1}.}} By \eqref{ut+} and \eqref{ut-},
\begin{equation*}
\jap{p}u(t,\pm)-e^{\pm it\jap{p}}\jap{p}u_{\Omega,\pm}=\mp\tfrac{1}{2i}\int_t^\infty e^{\pm i(t-s)\jap{p}}N(u(s),x,s)u(s)\,ds.
\end{equation*}
This, together with Lemma~\ref{Lem: Pprop}, gives
\begin{align*}
&\bigl\|F(|x|>t^\alpha)\,P^+_\pm F_{t^{-\beta}}(v(p))\bigl(\jap{p}u(t,\pm)-e^{\pm it\jap{p}}\jap{p}u_{\Omega,\pm}\bigr)\bigr\|\nonumber\\
&\leq \tfrac{1}{2}\int_t^\infty\bigl\|F(|x|>t^\alpha)P^+_\pm F_{t^{-\beta}}(v(p))e^{\pm i(t-s)\jap{p}}\langle x\rangle^{-\sigma}\bigr\|\,\bigl\|\langle x\rangle^\sigma N(u(s),x,s)u(s)\bigr\|\,ds\nonumber\\
&\leq \int_t^\infty\!\frac{CE_u}{\langle t^\alpha+|t-s|t^{-\beta}\rangle^\sigma}\,ds\leq \Bigl(\frac{C}{t^{\alpha\sigma-1}}+\frac{Ct^\beta}{t^{(1-\beta)(\sigma-1)}}\Bigr)E_u\to 0
\end{align*}
as $t\to\infty$, since $\alpha>1/\sigma$ and $\beta<1-1/\sigma$, where we used
\begin{equation*}
\int_t^\infty\!\frac{ds}{\langle t^\alpha+|t-s|t^{-\beta}\rangle^\sigma}\leq \int_t^{2t}\!\frac{ds}{\langle t^\alpha+|t-s|t^{-\beta}\rangle^\sigma}+\int_{2t}^\infty\!\frac{ds}{\langle t^\alpha+|t-s|t^{-\beta}\rangle^\sigma}\leq \frac{C}{t^{\alpha\sigma-1}}+\frac{Ct^\beta}{t^{(1-\beta)(\sigma-1)}}.
\end{equation*}

\noindent\fbox{\emph{Proof of \eqref{weak: eq3} for $\xi(t)=\jap{p}u(t,\pm)$.}} By Duhamel's formula applied to \eqref{ut+}--\eqref{ut-},
\eq\label{key: differ backward}
\jap{p}u(t,\pm)-e^{\pm it\jap{p}}\jap{p}u(0,\pm)=\pm\tfrac{1}{2i}\int_0^t e^{\pm i(t-s)\jap{p}}N(u,x,s)u(s)\,ds.
\eeq
By Lemma~\ref{Lem: Pprop} and Fubini's theorem, the integral part of $\jap{p}u(t,\pm)$ in \eqref{key: differ backward} satisfies
\begin{align}\label{est: weak: eq3 new}
&\bigl\|F(|x|>t^\alpha)P^-_\pm F_{t^{-\beta}}(v(p))\bigl(\jap{p}u(t,\pm)-e^{\pm it\jap{p}}\jap{p}u(0,\pm)\bigr)\bigr\|\nonumber\\
&\leq \tfrac{1}{2}\int_0^t\bigl\|F(|x|>t^\alpha)P^-_\pm F_{t^{-\beta}}(v(p))e^{\pm i(t-s)\jap{p}}\langle x\rangle^{-\sigma}\bigr\|\,\bigl\|\langle x\rangle^\sigma Nu(s)\bigr\|\,ds\nonumber\\
&\leq \int_0^t\!\frac{CE_u}{\langle t^\alpha+|t-s|t^{-\beta}\rangle^\sigma}\,ds\leq \frac{CE_u\,t}{t^{\alpha\sigma}}\to 0
\end{align}
as $t\to\infty$, since $\alpha>1/\sigma$. Combining \eqref{est: weak: eq3 new} (via \eqref{key: differ backward}) with Lemma~\ref{lem: low freq conv} yields \eqref{weak: eq3}.

\noindent\fbox{\emph{Proof of \eqref{weak: eq4} for $\xi(t)=\jap{p}u(t,\pm)$.}} To use Corollary~\ref{cor: small}, take $\beta=\tfrac{1}{2}$ and $\alpha>\tfrac{1}{2}$. By \eqref{key: differ backward} and Corollary~\ref{cor: small}, the integral part of $\jap{p}u(t,\pm)$ is controlled by an analogue of \eqref{est: weak: eq3 new} with $P^-_\pm$ removed:
\begin{align}\label{est: weak: eq5 int}
&\bigl\|F(|x|>t^\alpha)\,(1-F_{t^{-\beta}}(v(p)))\,\tfrac{1}{2i}\!\int_0^t e^{\pm i(t-s)\jap{p}}N(u,x,s)u(s)\,ds\bigr\|\nonumber\\
&\leq \tfrac{1}{2}\int_0^t\bigl\|F(|x|>t^\alpha)(1-F_{t^{-\beta}}(v(p)))e^{\pm i(t-s)\jap{p}}\langle x\rangle^{-\sigma}\bigr\|\,\bigl\|\langle x\rangle^\sigma Nu(s)\bigr\|\,ds\leq \frac{CE_u\,t}{t^{\alpha\sigma}}\to 0
\end{align}
as $t\to\infty$, since $\alpha>1/\sigma$. The free part $e^{\pm it\jap{p}}\jap{p}u(0,\pm)$ also satisfies, by Corollary~\ref{cor: small} (cf.~the integrand factor of \eqref{est: weak: eq5 int} with $s=0$),
\begin{align}\label{est: weak: eq5 free}
&\bigl\|F(|x|>t^\alpha)(1-F_{t^{-\beta}}(v(p)))e^{\pm it\jap{p}}\jap{p}u(0,\pm)\bigr\|\nonumber\\
&\leq \frac{C}{\langle t^\alpha+t^{1-\beta}\rangle^\sigma}\bigl\|\chi(|x|<t^{\alpha/2})\langle x\rangle^\sigma\jap{p}u(0,\pm)\bigr\|+\bigl\|\chi(|x|\geq t^{\alpha/2})\jap{p}u(0,\pm)\bigr\|\nonumber\\
&\leq \frac{C\,t^{\alpha\sigma/2}}{\langle t^\alpha+t^{1-\beta}\rangle^\sigma}\bigl\|\jap{p}u(0,\pm)\bigr\|+\bigl\|\chi(|x|\geq t^{\alpha/2})\jap{p}u(0,\pm)\bigr\|\to 0
\end{align}
as $t\to\infty$. Combining \eqref{est: weak: eq5 int}--\eqref{est: weak: eq5 free} with \eqref{key: differ backward} yields \eqref{weak: eq4} for $\xi(t)=\jap{p}u(t,\pm)$.
\end{proof}
\begin{proof}[Proof of Theorem~\ref{thm3}]
By Lemma~\ref{lem: weak err}, for every $\alpha>\max\{\tfrac{1}{2},\tfrac{1}{\sigma}\}$, as $t\to\infty$,
\eq\label{step1: thm3}
\bigl\|F(|x|>t^\alpha)\bigl(\jap{p}\,u(t,\pm)-P^+_\pm F_{t^{-\beta}}(v(p))\,e^{\pm it\jap{p}}\jap{p}\,u_{\Omega,\pm}\bigr)\bigr\|\to 0,
\eeq
and similarly with $\dot u(t,\pm)$ in place of $\jap{p}\,u(t,\pm)$ and $\dot u_{\Omega,\pm}$ in place of $\jap{p}\,u_{\Omega,\pm}$. Indeed, the decomposition \eqref{decomp: u sub} gives
\begin{align*}
\jap{p}\,u(t,\pm)-P^+_\pm F_{t^{-\beta}}(v(p))\,&e^{\pm it\jap{p}}\jap{p}\,u_{\Omega,\pm}=P^+_\pm F_{t^{-\beta}}(v(p))\bigl(\jap{p}\,u(t,\pm)-e^{\pm it\jap{p}}\jap{p}\,u_{\Omega,\pm}\bigr)\\
&\ +P^-_\pm F_{t^{-\beta}}(v(p))\jap{p}\,u(t,\pm)+(1-F_{t^{-\beta}}(v(p)))\jap{p}\,u(t,\pm),
\end{align*}
and applying $F(|x|>t^\alpha)$ and taking norms, the three summands tend to $0$ by \eqref{weak: eq1}, \eqref{weak: eq3} and \eqref{weak: eq4}, respectively.

Next, we claim that as $t\to\infty$,
\eq\label{step2: thm3}
\bigl\|P^+_\pm F_{t^{-\beta}}(v(p))\,e^{\pm it\jap{p}}\jap{p}\,u_{\Omega,\pm}-e^{\pm it\jap{p}}\jap{p}\,u_{\Omega,\pm}\bigr\|\to 0,
\eeq
and similarly for the $\dot u_{\Omega,\pm}$ counterpart. Indeed, using $P^+_\pm+P^-_\pm=I$, the LHS equals $\|-P^+_\pm(1-F_{t^{-\beta}}(v(p)))\,e^{\pm it\jap{p}}\jap{p}\,u_{\Omega,\pm}-P^-_\pm e^{\pm it\jap{p}}\jap{p}\,u_{\Omega,\pm}\|$ and is thus bounded by
\begin{equation*}
\|(1-F_{t^{-\beta}}(v(p)))\jap{p}\,u_{\Omega,\pm}\|+\|P^-_\pm e^{\pm it\jap{p}}\jap{p}\,u_{\Omega,\pm}\|.
\end{equation*}
The first term tends to $0$ by Lemma~\ref{lem: low freq conv}. For the second, set $R(t):=t^{(1-\beta)/2}$ and split $\jap{p}\,u_{\Omega,\pm}=\chi(|x|<R(t))\jap{p}\,u_{\Omega,\pm}+\chi(|x|\geq R(t))\jap{p}\,u_{\Omega,\pm}$ as well as $F_{t^{-\beta}}(v(p))+(1-F_{t^{-\beta}}(v(p)))=I$; then by Lemma~\ref{Lem: Pprop2} and Lemma~\ref{lem: low freq conv},
\begin{align*}
&\|P^-_\pm e^{\pm it\jap{p}}\jap{p}\,u_{\Omega,\pm}\|\\
&\leq \frac{C\,\langle R(t)\rangle^\sigma}{(t^{1-\beta})^\sigma}\|\jap{p}\,u_{\Omega,\pm}\|+\|\chi(|x|\geq R(t))\jap{p}\,u_{\Omega,\pm}\|+\|(1-F_{t^{-\beta}}(v(p)))\jap{p}\,u_{\Omega,\pm}\|\to 0,
\end{align*}
which yields \eqref{step2: thm3}.

Combining \eqref{step1: thm3} and \eqref{step2: thm3} by the triangle inequality gives, for every $\alpha>\max\{\tfrac{1}{2},\tfrac{1}{\sigma}\}$, as $t\to\infty$,
\eq\label{step3: thm3}
\bigl\|F(|x|>t^\alpha)\bigl(\jap{p}\,u(t,\pm)-e^{\pm it\jap{p}}\jap{p}\,u_{\Omega,\pm}\bigr)\bigr\|\to 0,
\eeq
and the analogous statement holds with $\dot u(t,\pm)$ and $\dot u_{\Omega,\pm}$ in place of $\jap{p}\,u(t,\pm)$ and $\jap{p}\,u_{\Omega,\pm}$. By the representation \eqref{uOmega: lc} together with the identity $\dot u_{\Omega,\pm}=\pm i\jap{p}\,u_{\Omega,\pm}$, the components of $U_0(t,0)\Omega_\alpha^*\vec u(0)$ read
\eq\label{Omega: prop}
\jap{p}\,u_{\Omega,1}(t)=\sum_{\pm}e^{\pm it\jap{p}}\jap{p}\,u_{\Omega,\pm},\qquad u_{\Omega,2}(t)=\sum_{\pm}e^{\pm it\jap{p}}\dot u_{\Omega,\pm}.
\eeq
Summing \eqref{step3: thm3} (and its $\dot u$-analogue) over the two signs, and using the decomposition $u(t)=u(t,+)+u(t,-)$, $\dot u(t)=\dot u(t,+)+\dot u(t,-)$ from \eqref{decomp: u total} together with \eqref{Omega: prop}, we obtain \eqref{weak: goal1} and \eqref{weak: goal2}. Defining $\vec u_{wlc}(t):=F(|x|\leq t^\alpha)\vec u(t)$ then completes the proof of Theorem~\ref{thm3}.
\end{proof}

\section{Examples}

In this section we illustrate the main results by two concrete examples. The first is a one-dimensional KG equation with a local interaction as in Example~\ref{ex: local} (Corollary~\ref{app1}, of Theorems~\ref{thm1} and~\ref{thm3}). The second is a higher-dimensional KG equation with a charge-transfer potential plus power-type nonlinearities as in Example~\ref{ex: charge} (Corollary~\ref{app3}, of Theorem~\ref{thm2}).

\subsection{One-dimensional KG with a local interaction}\label{subsec: ex1}

We consider, in space dimension $n=1$, the equation \eqref{KG} with $N$ as in Example~\ref{ex: local}, namely
\eq\label{KG0}
(\square+1)u=V(x,t)u+a(x)u^2+b(x)u^3,\qquad \vec u(0)\in\Hi,\quad (x,t)\in \mathbb{R}\times \mathbb{R}_{\geq 0},
\eeq
where $\langle x\rangle^\sigma V(x,t)\in L^\infty_{x,t}(\mathbb{R}\times\mathbb{R}_{\geq 0})$ and $\langle x\rangle^\sigma a(x),\langle x\rangle^\sigma b(x)\in L^\infty_x(\mathbb{R})$ for some $\sigma>1$, so that the interaction
\begin{equation*}
N(u,x,t)=V(x,t)+a(x)u+b(x)u^2
\end{equation*}
is local in the sense of Assumption~\ref{asp: 2}.

\begin{corollary}[of Theorems~\ref{thm1} and~\ref{thm3}]\label{app1}Let $\vec u(t)$ be a global solution of \eqref{KG0} satisfying Assumption~\ref{asp: global}, with $V$, $a$, $b$ as above. Then $\vec u(t)$ verifies Assumption~\ref{asp: 2}, and consequently Theorems~\ref{thm1} and~\ref{thm3} apply: the free channel wave operator $\Omega_\alpha^*\vec u(0)$ exists in $\Hi$, and as $t\to\infty$,
\begin{equation*}
\|\vec u(t)-U_0(t,0)\Omega_\alpha^*\vec u(0)-\vec u_{wlc}(t)\|_\Hi\to 0,
\end{equation*}
where $\vec u_{wlc}(t)$ satisfies the sub-ballistic estimate \eqref{conlusion: wlocal}.
\end{corollary}
\begin{proof}It suffices to verify Assumption~\ref{asp: 2}. By Assumption~\ref{asp: global} and the Sobolev embedding $H^1(\mathbb{R})\hookrightarrow L^\infty(\mathbb{R})$ (valid in dimension $n=1$),
\begin{equation*}
\|u(t)\|_{L^\infty_x}\lesssim \|u(t)\|_{H^1_x}\lesssim_E 1.
\end{equation*}
Hence
\begin{align*}
\|\langle x\rangle^\sigma N(u(t),x,t)u(t)\|&\leq \|\langle x\rangle^\sigma V\|_{L^\infty_{x,t}}\|u(t)\|+\|\langle x\rangle^\sigma a\|_{L^\infty_x}\|u(t)\|_{L^\infty_x}\|u(t)\|\\
&\quad+\|\langle x\rangle^\sigma b\|_{L^\infty_x}\|u(t)\|_{L^\infty_x}^2\|u(t)\|\\
&\lesssim_E 1,
\end{align*}
which is Assumption~\ref{asp: 2}. The conclusion now follows from Theorems~\ref{thm1} and~\ref{thm3}.
\end{proof}

\subsection{Higher-dimensional KG with charge-transfer plus power nonlinearities}\label{subsec: ex2}

We now turn to space dimension $n\geq 3$ and consider the equation \eqref{KG} with $N$ as in Example~\ref{ex: charge}, namely
\eq\label{N charge}
N(u,x,t)=V(x,t)+\sum_{j=1}^{K}\epsilon_j\lambda_j|u|^{p_j},\qquad \lambda_j>0,\ \epsilon_j\in\{+1,-1\},\ 1\leq p_j\leq \frac{n+2}{n-2},
\eeq
where $V(x,t)$ is of charge-transfer type,
\eq\label{V charge}
V(x,t)=\sum_{j=1}^M V_j(x-g_j(t)v_j,t),\quad V_j\in L^\infty_tL^2_x(\mathbb{R}^n\times\mathbb{R}_{\geq 0}),\ v_j\in \mathbb{R}^n,
\eeq
with real-valued $g_j(t)$.

\begin{corollary}[of Theorem~\ref{thm2}]\label{app3}Let $n\geq 3$ and let $\vec u(t)$ be a global solution of \eqref{KG} with $N$ as in \eqref{N charge}--\eqref{V charge}, satisfying Assumption~\ref{asp: global}. Then $\vec u(t)$ verifies Assumption~\ref{asp: 3}, and consequently Theorem~\ref{thm2} applies: the free channel wave operator $\Omega_\alpha^*\vec u(0)$ exists in $\Hi$.
\end{corollary}
\begin{proof}It suffices to verify Assumption~\ref{asp: 3}. By the Sobolev embedding $H^1\hookrightarrow L^{2n/(n-2)}$ on $\mathbb{R}^n$, together with Assumption~\ref{asp: global},
\begin{equation*}
\||u(t)|^{p_j}u(t)\|_{L^1_x}=\|u(t)\|_{L^{p_j+1}_x}^{p_j+1}\lesssim_E 1
\end{equation*}
for each $1\leq p_j\leq (n+2)/(n-2)$, since $p_j+1\in [2,\tfrac{2n}{n-2}]$ lies in the admissible Sobolev range. Moreover, by Cauchy--Schwarz and \eqref{V charge},
\begin{equation*}
\|V(\cdot,t)u(t)\|_{L^1_x}\leq \|V(\cdot,t)\|\|u(t)\|\leq \Big(\sum_{j=1}^M \|V_j\|_{L^\infty_tL^2_x}\Big)\sup_{s\geq 0}\|u(s)\|<\infty.
\end{equation*}
Combining these, $\sup_{t\geq 0}\|Nu(t)\|_{L^1_x}<\infty$, which is Assumption~\ref{asp: 3}.
\end{proof}

\section*{Acknowledgments}
X.W.\ acknowledges support from the ARC Australian Laureate Fellowship, grant FL220100072. A.S.\ acknowledges support from the National Science Foundation, grant DMS-2205931. Part of this work was carried out while X.W.\ was at Rutgers University.

\bibliographystyle{amsplain-unsrt}
\providecommand{\bysame}{\leavevmode\hbox to3em{\hrulefill}\thinspace}
\providecommand{\MR}{\relax\ifhmode\unskip\space\fi MR }
\providecommand{\MRhref}[2]{%
  \href{http://www.ams.org/mathscinet-getitem?mr=#1}{#2}
}
\providecommand{\href}[2]{#2}

\end{document}